\documentclass{imamat}
\usepackage{graphicx}
\usepackage{bm}
\usepackage{tikz}
\usetikzlibrary{decorations.pathmorphing}
\usepackage{caption,subcaption}
\usepackage{float}
\usepackage[leqno]{subeqnarray}
\usepackage{slashbox}
\newtheorem{thm}{Theorem}[section]
\newcommand{\Weber}{\operatorname{\mathit{W\kern-.30em e}}}
\newcommand{\Rey}{\operatorname{\mathit{R\kern-.10em e}}}

\usepackage{multirow}
\usepackage{enumerate}

\jno{xxx000}
\received{25 January 2019}
\begin{document}

\title{Point-actuated feedback control of multidimensional interfaces}
% Short title for running heads:
\shorttitle{Point-actuated feedback control of multidimensional interfaces}

\author{%
{\sc
Ruben J. Tomlin\thanks{Corresponding author. Email: rjt111@ic.ac.uk}
and
Susana N. Gomes\thanks{Email: susana.gomes@warwick.ac.uk}}\\[2pt]
Department of Mathematics, Imperial College London
}
% Short list of authors for running heads:
\shortauthorlist{R. J. Tomlin and S. N. Gomes}

\maketitle

\begin{abstract}
% Body of abstract:
{We consider the application of feedback control strategies with point actuators to stabilise desired interface shapes. We take a multidimensional Kuramoto--Sivashinsky equation as a test case; this equation arises in the study of thin liquid films, exhibiting a wide range of dynamics in different parameter regimes, including unbounded growth and full spatiotemporal chaos. In the case of limited observability, we utilise a proportional control strategy where forcing at a point depends only on the local observation. We find that point-actuated controls may inhibit unbounded growth of a solution, if they are sufficient in number and in strength, and can exponentially stabilise the desired state. We investigate actuator arrangements, and find that the equidistant case is optimal, with heavy penalties for poorly arranged actuators. We additionally consider the problem of synchronising two chaotic solutions using proportional controls. In the case when the full interface is observable, we construct feedback gain matrices using the linearised dynamics. Such controls improve on the proportional case, and are applied to stabilise non-trivial steady and travelling wave solutions.}
% Keywords:
{feedback control, interfacial dynamics, point actuators, proportional control, multidimensional Kuramoto--Sivashinsky equation, thin films.}
\end{abstract}

\section{Introduction}

The study of evolving interfaces is at the core of many areas of applied mathematics, ranging from growth processes in mathematical biology \citep{eden1961two} and chemistry \citep{kobayashi1993modeling}, to the evolution of liquid--air interfaces in fluid dynamics \citep{Siv3}, flame front propagation in combustion theory \citep{Siv1,Siv2,Siv4}, or even defector/cooperator problems in game theory \citep{szolnoki2018evolutionary}. Starting with complicated multiphase systems comprising numerous coupled equations, it may be possible to employ modelling techniques to isolate the evolution of the interface alone; this is particularly desirable when information about the bulk dynamics (away from the interface) is not of interest. It is often challenging to extract the interfacial dynamics while still retaining all the desired physical effects, and in many cases it is found that the obtained low-dimensional models only replicate the true dynamics well in restricted parameter regimes. However, this is balanced by the relative simplicity of the interface evolution equations along with a large decrease in computational complexity for numerical simulations.

It may be useful to control the interfacial dynamics in order to optimise a process. For example, cooling and coating processes arise in microfluidic applications where a thin liquid film flows over a substrate. In the former case, waviness of the liquid interface is desirable as it improves heat and mass transfer \citep{LYU19911451,Miyara1999,serifi2004transient}, whereas in the latter case, a flat interface is needed. For a thin film flow, controls may take the form of air/liquid actuators, electric/magnetic fields, surfactants or substrate coating/topography. Controls may also be introduced for crystal growth processes where the rate of growth can be modified with heat sources \citep{kokh2005application}. Through the modelling procedure described above, controls at the level of the full physical system are recast as controls acting on the interface alone -- boundary controls manifest themselves as distributed (internal) controls acting on the interface.

For most real-world problems, interfaces are described in terms of two spatial variables. For problems where variations in one direction are negligible, such as the growth of a flat crystal, simplification of the interface problem to one spatial dimension may be viable, and it is important that the controls utilised preserve this property. However, it may be the case that such a simplification overlooks important instabilities or mechanisms which are only observed from the full three-dimensional (3D) formulation of the original multiphase problem, e.g. Rayleigh--Taylor instabilities or electrostatically induced instabilities in liquid films \citep{tomlin_papageorgiou_pavliotis_2017,tomlin2019optimal}.

This paper investigates two feedback control strategies for multidimensional interfaces using a 2D Kuramoto--Sivashinsky equation (KSE) as a test case. Controls are applied using a finite set of point-actuators. Such actuators are one of the most physically realisable, with localised forcing applied at specified nodes, injecting or extracting mass from the bulk flow which accordingly forces the interface. In the case of the 2D KSE under consideration here, which models the interface of a thin film flow over a flat substrate, such controls arise via same-fluid blowing and suction at the substrate surface. Throughout this work, the point-actuated controls are manifested mathematically as Dirac delta functions; smoothed alternatives have been utilised throughout the literature. Many authors have considered the use of point actuators for fluid interfaces in one spatial dimension \citep{christofides1998feedback,armaou2000feedback,lunasin2017finite,gomes2016stabilizing}.

The system under consideration in this work is an extension to two spatial dimensions of the 1D KSE,
\begin{equation} \label{1dintroks} 
\eta_t + \eta \eta_x + \eta_{xx}  + \eta_{xxxx}  =  0,
\end{equation}
which is the paradigmatic model for the class of active-dissipative nonlinear evolution equations. Usually, \eqref{1dintroks} is supplemented with periodic boundary conditions on the interval $[0,L]$. As $L$ is increased beyond $2\pi$ (at which point the first Fourier mode destabilises), the dynamics cascades to full spatiotemporal chaos through steady and travelling wave, time-periodic, and quasi-periodic attractors. In this paper, we consider control strategies for the KSE in two space dimensions,
\begin{equation} \label{controlled2dks} 
\eta_t + \eta \eta_x + (1- \kappa) \eta_{xx}  -  \kappa  \eta_{yy} + \Delta^2 \eta  =  \zeta,
\end{equation}
where $\zeta$ is the control. This equation may be derived to describe the weakly nonlinear evolution of small-amplitude, long-wave perturbations of gravity-driven thin liquid films on flat substrates; the surface $\eta(x,y,t)$ represents a perturbation of the flat film solution. We supplement \eqref{controlled2dks} with periodic boundary conditions on the rectangular domain $Q = [0,L_1]\times [0,L_2]$. Different dynamical regimes are found by varying $\kappa$. Omitting chaotic dynamics or unbounded growth, the 2D KSE \eqref{controlled2dks}, albeit a deterministic equation, provides a natural and challenging test case for the control strategies considered in this work. Efficient and convergent numerical schemes allow us to study the control of \eqref{controlled2dks} on large domains with many unstable modes and solutions exhibiting full spatiotemporal chaos; the majority of existing numerical studies consider parameter regimes for model problems not far from the onset of instability (one or two unstable modes).

The two feedback control (closed-loop) strategies we consider are proportional control and feedback control with full state observations (shortened to``full feedback control"). These strategies are polar opposites in terms of the assumed observability of the interface and knowledge of the governing dynamics. Unsurprisingly, we find that more information (observations/knowledge of governing equation) results in a much improved control performance, but both methods are effective. Note that the open-loop optimal control problem for \eqref{controlled2dks} was considered in \cite{tomlin2019optimal}.

Proportional controls are the most simplistic and physically realisable form of feedback control. Each actuator is paired with an observer, and the forcing applied by that actuator is proportional to the difference between the observation of the interface and the chosen desired state (e.g. a travelling wave solution). In this study, each observer is co-located with an actuator for simplicity, and these are paired for proportional control. We note that upstream (phase-shifted) observers were found to improve the control of thin film models in \cite{thompson2016stabilising}. No information of the dynamical system is required or even beneficial since the actuation at a particular point does not utilise observations from other spatial locations. We perform a number of numerical experiments using proportional controls, investigating different actuator arrangements and the control of exponentially growing or chaotic interfaces. Furthermore, we investigate the use of proportional controls to synchronise two chaotic solutions, having possible applications in communications -- see \cite{pecora2015synchronization} and the references therein.

Full feedback entails the more advanced closed-loop control strategy which involves observation of the full interface and assumes knowledge of the governing system. The linearisation of the dynamics is used to create a function (the feedback gain matrix) which maps the observation of the full state at an instant in time to the controls required to obtain linear stability of the desired state. We present two methodologies for full feedback control to non-trivial interface shapes based on work on the 1D KSE~\eqref{1dintroks} by \cite{doi:10.1137/140993417} and \cite{gomes2016stabilizing}. The former ensures exponential stabilisation through a rigorous analytical result, whereas the latter is much more feasible numerically if the desired state is non-trivial. The capabilities of full feedback control in this multidimensional setting are tested with comparisons against the proportional control results.

There is a long list of extensions and hybridisations of these control strategies which we do not consider in the current work, such as dynamical observers, time-delayed/phase-shifted observers, or feedback strategies where controls depend on different subsets of the observers. However, the present study considers methods which are readily extendable to more complicated systems and experiments.

The current paper is organised as follows: Section \ref{PhysicalModel} introduces the control problem for \eqref{controlled2dks} with a brief discussion of its relevance to fluid dynamics. We provide the analytical setting of the problem, and continue to discuss the numerical methods utilised and the various actuator arrangements considered. Sections \ref{SecPropCont} and \ref{SecFeedbackcontrol} contain the studies of the proportional and full feedback control strategies, respectively. The concluding remarks are given in Section \ref{ConcSec}.

\section{Multidimensional Kuramoto--Sivashinsky equation with point-actuated controls\label{PhysicalModel}}

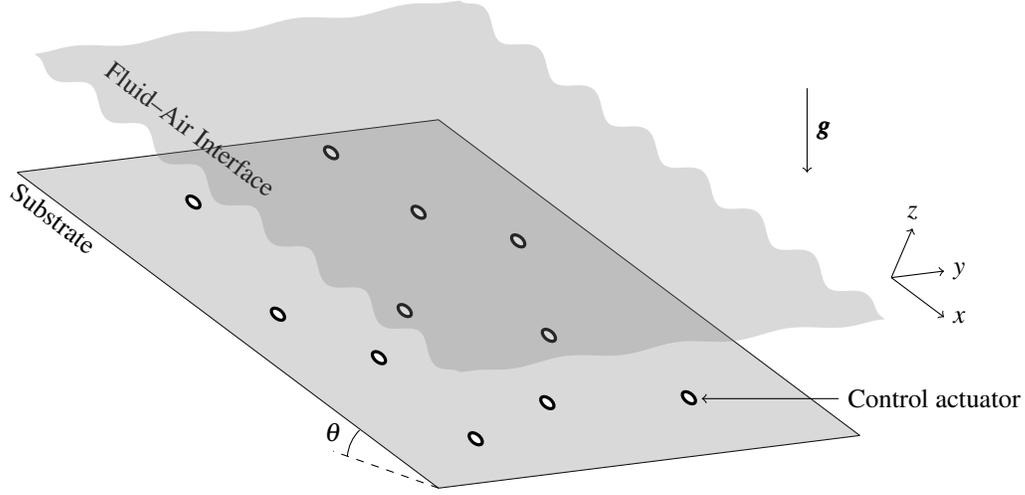
\begin{figure}
\begin{tikzpicture}[scale=1.40]

    \draw (5,0) coordinate (A1) -- (1,3) coordinate (A2);
    \draw (9,0.5) coordinate (A3) -- (5,3.5) coordinate (A4);
    \draw (A1) -- (A3); \draw (A2) -- (A4);
    \fill[gray,opacity=0.3] (A1) -- (A2) -- (A4) -- (A3) -- cycle;
    \draw[dashed] (A1) -- (4,0.35) node[above] {$\theta$};
    \draw (4.25,0.55) arc (135:175:0.4);

  \draw[-] (1.8, 4.05) node[right, rotate=-37] {Fluid--Air Interface};%{$z=h(x,y,t)$};
    
  \draw[-] (0.9, 2.9) node[right, rotate=-37] {Substrate};%{$z=0$};

    \draw[dashed, opacity=0.0] (5.2,1.466-0.35) coordinate (B1) -- (1.2,4.466-0.35) coordinate (B2);
    \draw[dashed, opacity=0.0] (9.2,1.966-0.35) coordinate (B3) -- (5.2,4.966-0.35) coordinate (B4);
    \draw[dashed, opacity=0.0] (B1) -- (B3); 
    \draw[dashed, opacity=0.0] (B2) -- (B4);

    \draw[->] (9.3, 2) -- (9.8,2.0625) node[right] {$y$};
    \draw[->] (9.3, 2) -- (9.8, 1.625) node[right] {$x$};
    \draw[->] (9.3, 2) -- (9.5, 2.466) node[above] {$z$};

     \draw[rotate=140,line width=1.2pt,fill=white] (-0.3,-3.8) ellipse (2.1pt and 1.2pt);
       
     \draw[rotate=140,line width=1.2pt,fill=white] (-2.6,-3.8) ellipse (2.1pt and 1.2pt);
     
     \draw[rotate=140,line width=1.2pt,fill=white] (-2,-5.1) ellipse (2.1pt and 1.2pt);
       
     \draw[rotate=140,line width=1.2pt,fill=white] (-5.1,-5.4) ellipse (2.1pt and 1.2pt);
     
      \draw[rotate=140,line width=1.2pt,fill=white] (-1.6,-3.5) ellipse (2.1pt and 1.2pt);
       
     \draw[rotate=140,line width=1.2pt,fill=white] (-1,-5) ellipse (2.1pt and 1.2pt);
     
      \draw[rotate=140,line width=1.2pt,fill=white] (-4.1,-4.5) ellipse (2.1pt and 1.2pt);
      
      \draw[rotate=140,line width=1.2pt,fill=white] (-2.9,-5.5) ellipse (2.1pt and 1.2pt);
   
  \draw[rotate=140,line width=1.2pt,fill=white] (-3.8,-3.8) ellipse (2.1pt and 1.2pt);
    
  \draw[rotate=140,line width=1.2pt,fill=white] (-3.7,-5) ellipse (2.1pt and 1.2pt);
   
  \draw[rotate=140,line width=1.2pt,fill=white] (-2.5,-4.3) ellipse (2.1pt and 1.2pt);

%         \def\mypath{ (B1) [decorate, decoration={snake, amplitude=0.07cm,segment length=0.911cm}] -- (B2) -- (B4) -- (B4) -- (B3) -- (B1) }
%         
%             \fill[gray,opacity=0.2] \mypath;
         
   \filldraw[gray,opacity=0.3, very thick]                    (B1)
  decorate [decoration={snake,amplitude=0.12cm,segment length=0.911cm}]             { -- (B2) }
  decorate [decoration={snake,amplitude=0.04cm,segment length=1.5cm}] 		{ -- (B4) }
  decorate [decoration={snake,amplitude=0.12cm,segment length=0.911cm}]           { -- (B3) }
   decorate [decoration={snake,amplitude=0.04cm,segment length=1.5cm}]           { -- (B1) }
      ;
    
 \draw[->] (8.5, 3.8) -- (8.5, 3.0);
     
      \draw[-] (8.5, 3.4) node[right] {$\bm{g}$};
      
      %Control
      
\draw[->] (8.8,0.85) -- (7.5,0.85);
     
      \draw[-] (8.8,0.85) node[right] {Control actuator};

\end{tikzpicture}\caption{Schematic for point-actuated control of an overlying thin liquid film.} \label{Setupdiagram}
\end{figure}

We consider the feedback control problem for the 2D KSE \eqref{controlled2dks}. This equation may be derived, with the addition of an advective term of the form $\chi \eta_x$, i.e.
\begin{equation} \label{controlled2dksadv} 
\eta_t + \chi \eta_x + \eta \eta_x + (1- \kappa) \eta_{xx}  -  \kappa  \eta_{yy} + \Delta^2 \eta  =  \zeta,
\end{equation}
in the context of thin film flows on inclined flat substrates -- see \cite{tomlin2019optimal} for example. The schematic in Figure \ref{Setupdiagram} shows the set-up for an overlying thin film flow with same-fluid blowing and suction controls at the substrate surface. The position in the $z$-coordinate of the fluid interface is a function of the streamwise and transverse spatial variables, $x$ and $y$, respectively, and time $t$. In the setting of Figure \ref{Setupdiagram}, $\eta$ represents a scaled interfacial perturbation, and $\zeta$ represents the forcing on the interface due to actuation. The parameter $\kappa$ encodes the angle of the substrate to the horizontal: For $\kappa > 0$ we have overlying film flows, a vertical film flow for $\kappa = 0$, and hanging flows when $\kappa < 0$. The value of $\kappa = 1$ corresponds to taking the critical Reynolds number, ${\Rey}_{\textrm{c}}$, for an overlying flow. As $\kappa$ is decreased from $1$, the flat film solution first becomes unstable to long waves (${\Rey} > {\Rey}_{\textrm{c}} $); for $\kappa > 1$ we have subcritical Reynolds number flows (${\Rey}< {\Rey}_{\textrm{c}} $), where all initial conditions decay to zero uniformly in the absence of controls. The aforementioned advection parameter $\chi$ measures the speed of waves travelling downstream relative to the ``lab" frame of reference. In the thin film context, $\chi$ is large. If the set of actuators is invariant to shifts in the streamwise direction (which is not the case in the current study), it is viable to employ the Galilean transformation $x \rightarrow x + \chi t$ (substitute $x = \overline{x} + \chi t$ into \eqref{controlled2dksadv} and drop the bar) to remove the advective term, resulting in equation \eqref{controlled2dks}. The advection term provides no additional complexity to the problem, and the fluid dynamicist's perspective alone merits the study of its effect on the controlled dynamics, since the results for \eqref{controlled2dks} are not immediately applicable to the thin film scenario. We found that increasing $\chi$ improved the controlled dynamics (i.e. increased decay rates), thus the case of $\chi = 0$ is seemingly the most challenging. We thus continue without the linear advective term; it is noteworthy that previous studies of the 1D KSE with point-actuated controls also ignored advection effects.

We supplement \eqref{controlled2dks} with periodic boundary conditions on the rectangle $Q = [0,L_1]\times[0,L_2]$, in which case the spectrum of solutions is restricted to wavenumbers $\bm{\tilde{k}} = (\tilde{k}_1,\tilde{k}_2)$ where
\begin{equation}\label{ktildedefn}\tilde{k}_1 = \frac{2\pi k_1}{L_1},\qquad \tilde{k}_2 = \frac{2\pi k_2}{L_2},\end{equation}
for $\bm{k} \in \mathbb{Z}^2$. We may write $\eta$ and $\zeta$ in terms of Fourier series as
\begin{equation}\label{fourierseriesetazeta1}\eta = \sum_{\bm{k}\in\mathbb{Z}^2} \eta_{\bm{k}} e^{i \bm{\tilde{k}} \bm{\cdot} \bm{x}}, \qquad \zeta = \sum_{\bm{k}\in\mathbb{Z}^2} \zeta_{\bm{k}} e^{i \bm{\tilde{k}} \bm{\cdot} \bm{x}}.\end{equation}
Equation \eqref{controlled2dks} may thus be replaced by an equivalent system of ODEs for the Fourier coefficients $\eta_{\bm{k}}$ and $\zeta_{\bm{k}}$,
\begin{equation}\label{fouriercoeffsfull1}\frac{\mathrm{d}}{\mathrm{d}t}\eta_{\bm{k}} + \frac{i\tilde{k}_1}{2} \sum_{\bm{m}\in\mathbb{Z}^2} \eta_{\bm{k}-\bm{m}} \eta_{\bm{m}} = \left[ (1-\kappa) \tilde{k}_1^2 - \kappa \tilde{k}_2^2 - |\bm{\tilde{k}}|^4  \right] \eta_{\bm{k}} + \zeta_{\bm{k}},\end{equation}
where $\eta_{-\bm{k}}$ and $\zeta_{-\bm{k}}$ are the complex conjugates of $\eta_{\bm{k}}$ and $\zeta_{\bm{k}}$, respectively, since both the solution and control are real-valued. The effect of the parameter $\kappa$ can be seen more clearly from this ODE system, and the different dynamical regimes are given in the following table:
\begin{table}[H]
\centering
\begin{tabular}{| c  | c | c | c | c |}
  \hline			
  $\kappa < 0$ & $\kappa = 0$ & $0 < \kappa < 1$ & $\kappa = 1$ & $1 < \kappa$ \\
  \hline
  hanging films & vertical film &    \multicolumn{3}{c|}{$\;\;\;\;$overlying films}  \\ \cline{3-5} \hline
     \multicolumn{3}{|c|}{$\;\;\quad\qquad{\Rey} > {\Rey}_{\textrm{c}}$} & ${\Rey} = {\Rey}_{\textrm{c}}$ & ${\Rey} < {\Rey}_{\textrm{c}}$ \\
  \hline
    unbounded growth   &  \multicolumn{2}{c|}{bounded non-trivial dynamics}     & \multicolumn{2}{c|}{flat solution is stable} \\
  \hline  
\end{tabular}
\caption{Regimes of dynamics for (\ref{controlled2dks},\ref{fouriercoeffsfull1}) in terms of $\kappa$.}\label{tab:kapparanges}
\end{table}
%\begin{table}[H]
%\centering
%\begin{tabular}{| c  | c | c | c | c |}
%  \hline			
%  $\kappa < 0$ & $\kappa = 0$ & $0 < \kappa < 1$ & $\kappa = 1$ & $1 < \kappa$ \\
%  \hline
%  hanging films & vertical film &    \multicolumn{3}{c|}{overlying films}  \\ \cline{3-5}
%   & &  ${\Rey} > {\Rey}_{\textrm{c}}$ & ${\Rey} = {\Rey}_{\textrm{c}}$ & ${\Rey} < {\Rey}_{\textrm{c}}$ \\
%  \hline
%    unbounded growth   &  \multicolumn{2}{c|}{bounded non-trivial dynamics}     & \multicolumn{2}{c|}{flat solution is stable} \\
%  \hline  
%\end{tabular}
%\caption{Regimes of dynamics for (\ref{controlled2dks},\ref{fouriercoeffsfull1}) in terms of $\kappa$.}\label{tab:kapparanges}
%\end{table}
\noindent The unbounded growth for $\kappa < 0$ is due to a linear (Rayleigh--Taylor) instability in transverse modes (assuming $L_2$ is sufficiently large) which is not saturated by the nonlinearity -- this can be seen from the ODE system \eqref{fouriercoeffsfull1} where the nonlinear (summation) term has no contribution if $\tilde{k}_1 = 0$. Further details of the dynamical regimes of (\ref{controlled2dks},\ref{fouriercoeffsfull1}) can be found in \cite{tomlin2019optimal}; in particular, the vertical falling film case has been considered extensively both in analytical and numerical studies \citep{Nepo1,Nepo2,pinto1,pinto2,akrivislinearly,tkp2018}. We focus our attention on stabilising the dynamics in the unstable regimes, $\kappa < 1$. In contrast, for $\kappa \geq 1$, the relevant control problem is in destabilising the flat interface. It is evident from \eqref{fouriercoeffsfull1} that the solution mean is preserved by the dynamics if the controls are zero mean, otherwise $\eta_{\bm{0}}$ will drift.

We denote the desired state by $\overline{\eta}$; this will usually be the trivial zero solution, but we also consider travelling waves and fully chaotic solutions of the uncontrolled system. Although not always necessary, we assume throughout that $\overline{\eta}$ is an exact solution (stable or unstable) of the uncontrolled equation; convergence to non-solutions usually requires a dense set of actuators. Furthermore, the full feedback control methodology requires this assumption. We consider point-actuated controls, with actuators located at $\{\bm{x}_j\}_{j =1}^{N_{\textrm{ctrl}}} \subset Q$. The actuator function and time-dependent control corresponding to the location $\bm{x}_j$ are denoted by $b^j(\bm{x})$ and $\phi^j(t)$, respectively. In the thin film setting, these controls correspond to blowing or suction applied through holes in the substrate surface as depicted in Figure \ref{Setupdiagram}. In general, the control may be expressed as
\begin{equation}\label{PAcontrolform1}\zeta(\bm{x},t) = \sum_{j =1}^{N_{\textrm{ctrl}}} \phi^{j}(t) b^{j}(\bm{x}) , \quad\textrm{where} \quad b^{j} (\bm{x})= \sum_{\bm{k}\in\mathbb{Z}^2} b^{j}_{\bm{k}} e^{i \bm{\tilde{k}} \bm{\cdot} \bm{x}}.\end{equation}
The Fourier coefficients of $\zeta$ may thus be expressed as a linear combination of the $b^{j}_{\bm{k}}$ depending on the controls $\phi^{j}$. Many studies of point-actuated controls employ smoothed actuators which are centred at the given actuator locations -- for example Gaussians, or functions which can be obtained from rescalings and translations of
 $ \exp (( \cos x - 1) w^{-2})$
%\begin{equation}b^{j}(x) = \exp \left( \frac{\cos(x) - 1}{w^2}\right)\end{equation}}
as used by \cite{thompson2016stabilising} in their study of long-wave (Benney, weighted-residual) models related to the KSE in 1D. Note that the latter approximation converges to a $2\pi$-periodic extension of the usual Dirac delta $\delta(x)$ as $w \rightarrow 0$. We do not make such an analytic approximation, and consider the Dirac delta actuators in two space dimensions, $b^{j}(\bm{x}) = \delta(\bm{x} - \bm{x}_{j})$ for $j =1, 2, \ldots, N_{\textrm{ctrl}}$, with Fourier coefficients
\begin{equation}\label{deltafseriesdefn1}b^{j}_{\bm{k}} = \frac{1}{|Q|} e^{-i \bm{\tilde{k}} \bm{\cdot} \bm{x}_{j}}.\end{equation}
The division by $|Q| = L_1L_2$ in this expression ensures that the corresponding distribution is dimensionless so that
\begin{equation} \int_Q b^{j}(\bm{x}) v(\bm{x}) \; \mathrm{d}\bm{x} = v(\bm{x}_j)
\end{equation}
for $v \in H^2$, the Sobolev space of periodic functions with both first and second spatial derivatives in $L^2$ (note that $\delta$ is in the dual space $H^{-2}$ and $H^2 \subset C^0$ in 2D). It follows (by taking $v=1$) that the spatial integral of $b^{j}$ over $Q$ is well-defined and is unity; additionally, any truncations of $b^j$ which include the zero mode have unit spatial integral. For numerical experiments, we truncate the Fourier series of the control at the same refinement as for the solution. In this way, the two limits of improving the resolutions of the solution and control are taken together. The forcing is highly singular, and requires many Fourier modes for good spatial convergence.

We define the spatial $L^2$-inner product and corresponding norm as
\begin{equation}\label{L2norms}\left\langle v , w \right\rangle_{L_2} =  \frac{1}{|Q|} \int_Q v w \; \mathrm{d}\bm{x} = \sum_{ \bm{k} \in\mathbb{Z}^2 } v_{\bm{k}}  w_{-\bm{k}}, \qquad \| v \|_{L^2}^2 = \left\langle v , v \right\rangle_{L_2} = \sum_{ \bm{k} \in\mathbb{Z}^2 } |v_{\bm{k}}|^2,
\end{equation}
where $v_{\bm{k}}$ and $w_{\bm{k}}$ are the Fourier coefficients of $v$ and $w$, respectively, as in \eqref{fourierseriesetazeta1}. The particular scaling in \eqref{L2norms} gives meaning to the $L^2$-norm as a measure of interfacial energy density. The concerns of the current work are primarily numerical, however we make some brief remarks on the analytical aspects of the problem. In 2D, we have $\delta \in H^{-2}$, and thus for $\phi^j\in L^2(0,T)$ we conclude that $\zeta \in L^2(0,T;H^{-2})$, i.e. the $H^{-2}$-norm of $\zeta$ is $L^2$-in-time. This is the minimal regularity of forcing needed for existence and uniqueness of solutions to \eqref{controlled2dks} in $L^2(0,T;H^2) \cap C^0([0,T];L^2)$, assuming $\eta_0 \in L^2$. This is due to standard results for parabolic equations -- see Ch. III \S 3 of \cite{temam2001navier} or Ch. 9.4 of \cite{robinson2001infinite} for the corresponding results for the Navier--Stokes equations, with extension to KS-type equations following \cite{temam1997infinite}. The authors have also considered the possibility of obtaining analytical estimates for the number of controls and control strength sufficient for exponential stabilisation of \eqref{controlled2dks} in the proportional control case, as done for similar problems in 1D by \cite{azouani2014feedback,lunasin2017finite}. We observed for the 1D KSE \eqref{1dintroks} that the analytical result appears far from optimal. The extension to multiple spatial dimensions is not trivial. This discussion of analytical aspects is also fully relevant for the 2D KSE with advection.

\subsection{Numerical Methods and Data Analysis\label{Numerical Methodssubsec}} 

For our numerical study of \eqref{controlled2dks} on $Q$-periodic domains, we utilise backwards differentiation formula (BDF) methods for the time discretisation and spectral methods in space. The BDFs belong to the family of implicit--explicit methods constructed by \cite{akrivis2004linearly} for a class of nonlinear parabolic equations -- see the appendix of \cite{akrivis2009linearly} for the first- to sixth-order schemes. They considered evolution equations of the form 
\begin{equation}\label{numericsparab1}\eta_t + \mathcal{A}\eta = \mathcal{B}(\eta),\end{equation}
where $\mathcal{A}$ is a positive definite, self-adjoint linear operator, and $\mathcal{B}$ is a nonlinear operator which satisfies a local Lipschitz condition. It was shown that these numerical schemes are efficient, convergent, and unconditionally stable. The applicability of these schemes for \eqref{controlled2dks} without controls was shown in \cite{akrivis2011linearly} and a convergence study was performed in \cite{akrivislinearly} for the choice of $\kappa = 0$. It was observed that the BDF schemes of order three to six (which are not unconditionally stable) achieved convergence to machine accuracy as soon as the time-step was small enough for the stability of the scheme. The BDFs were also utilised in \cite{tomlin_papageorgiou_pavliotis_2017} for a non-local problem. These schemes have been employed for both the 1D and 2D optimal control problems for the KSE in \cite{gomes2016stabilizing} and \cite{tomlin2019optimal}, respectively. We predominantly utilised the fourth-order BDF scheme, and performed tests with the other schemes for validation. For us, with the addition of the forcing, the operators in \eqref{numericsparab1} are defined as
\begin{subeqnarray}\label{AandBdefns1}  
\gdef\thesubequation{\theequation \text{a,b}}
\mathcal{A}\eta = (1- \kappa) \eta_{xx}  -  \kappa  \eta_{yy} + \Delta^2 \eta + c\eta , \qquad \mathcal{B}(\eta,\zeta)  =  - \eta \eta_x + \zeta +  c\eta, 
\end{subeqnarray}
where $c$ is chosen to ensure that $\mathcal{A}$ is positive definite. The forcing $\zeta$, being a summation of Dirac delta functions, is very singular, although $\mathcal{B}$ still satisfies the required Lipschitz condition if the controls $\phi^j$ themselves are Lipschitz in the state $\eta$. This can be seen from the following calculation where we assume that the controls consist of one point actuator at $\bm{x} = \bm{0}$ with $\zeta(\eta) = \phi(t;\eta)\delta(\bm{x})$: for $v \in H^2$ we have
\begin{equation}\left\langle \zeta(\eta_1) - \zeta(\eta_2) , v \right\rangle_{L^2} = \frac{1}{|Q|} \int_Q [\phi(t;\eta_1) - \phi(t;\eta_2)]\delta(\bm{x}) v(\bm{x}) \; \mathrm{d}\bm{x}  = \frac{v(\bm{0})}{|Q|} [\phi(t;\eta_1) - \phi(t;\eta_2)],\end{equation}
where $v(\bm{0})$ is bounded by the $H^2$-norm of $v$. Thus, it only remains to check that the controls $\phi^j$ satisfy
\begin{equation} |\phi^j(t;\eta_1) - \phi^j(t;\eta_2)| \leq \mu \| \eta_1 - \eta_2 \|_{L^2}\end{equation}
for the control schemes, where $\mu$ is a Lipschitz constant. This is not true in general for the point observation cases, but is trivial for the full feedback control case since the methodology follows a Fourier series framework. However, since our controlled solutions remain sufficiently regular, this Lipschitz bound is not an issue; a convergence study follows in the next section.

We discretise the spatial domain $Q$ with $2M$ equidistant points in the streamwise $x$-direction, and $2N$ equidistant points in the transverse $y$-direction, producing a grid of $2M\times 2N$ spatial points, and a corresponding frequency truncation in Fourier space that resolves modes with wavenumbers $|k_1| \leq M-1$ and $|k_2| \leq N-1$. The BDF method is then applied to the 2D Fast Fourier Transform (FFT) of the discretised interface, and the nonlinearity is calculated using the 2D FFT of $\eta^2$ (for this, we note that $(\eta^2)_x/2 = \eta\eta_x$). With \eqref{PAcontrolform1} and \eqref{deltafseriesdefn1}, it is clear that actuators do not need to be centered at computational grid points; additionally, having computed the Fourier coefficients of the solution at a given time, observations at any point location in $Q$ may be obtained to spectral accuracy using \eqref{fourierseriesetazeta1}. However, taking observations at grid points reduces the computational cost of each time-step. In many parts we take random initial conditions for our numerical simulations; for these we use
\begin{equation}\eta_0(\bm{x}) =  \sum_{| \bm{k} |_{\infty}=1}^{20}  
 a_{\bm{k}} \cos(\bm{\tilde{k}} \bm{\cdot} \bm{x}  ) + b_{\bm{k}} \sin(\bm{\tilde{k}} \bm{\cdot} \bm{x} ),
\end{equation}
where the coefficients $a_{\bm{k}}$ and $b_{\bm{k}}$ are random numbers from the interval $(-0.05,0.05)$. All initial conditions we consider have zero spatial mean (the mean is a conserved quantity for the uncontrolled system).

We track the time-dependent costs
\begin{equation}\label{costs1_2_define}\mathcal{C}_1(t) = \|  \eta(\cdot,t) - \overline{\eta}(\cdot,t)  \|_{L^2}, \qquad  \mathcal{C}_2(t) = \frac{1}{|Q|} \sum_{j=1}^{N_{\textrm{ctrl}}}  |\phi^j|,\end{equation}
where $\mathcal{C}_1$ measures the cost of the solution deviation from the desired state, and $\mathcal{C}_2$ measures the cost of the controls. Both costs are spatially dimensionless which is appropriate for our study on spatially periodic domains -- the costs over one period are the same as the costs over any number of periods considered together. The costs $\mathcal{C}_1$ and $\mathcal{C}_2$ are not necessarily equivalent (in the analytical sense), however they yield the same behaviour in our numerical simulations as will be seen below.

\subsection{Actuator arrangements/grids}

\begin{figure}%[H]
\centering
\begin{subfigure}{2.8in}
\caption{Equidistant.} 
\includegraphics[width=2.8in]{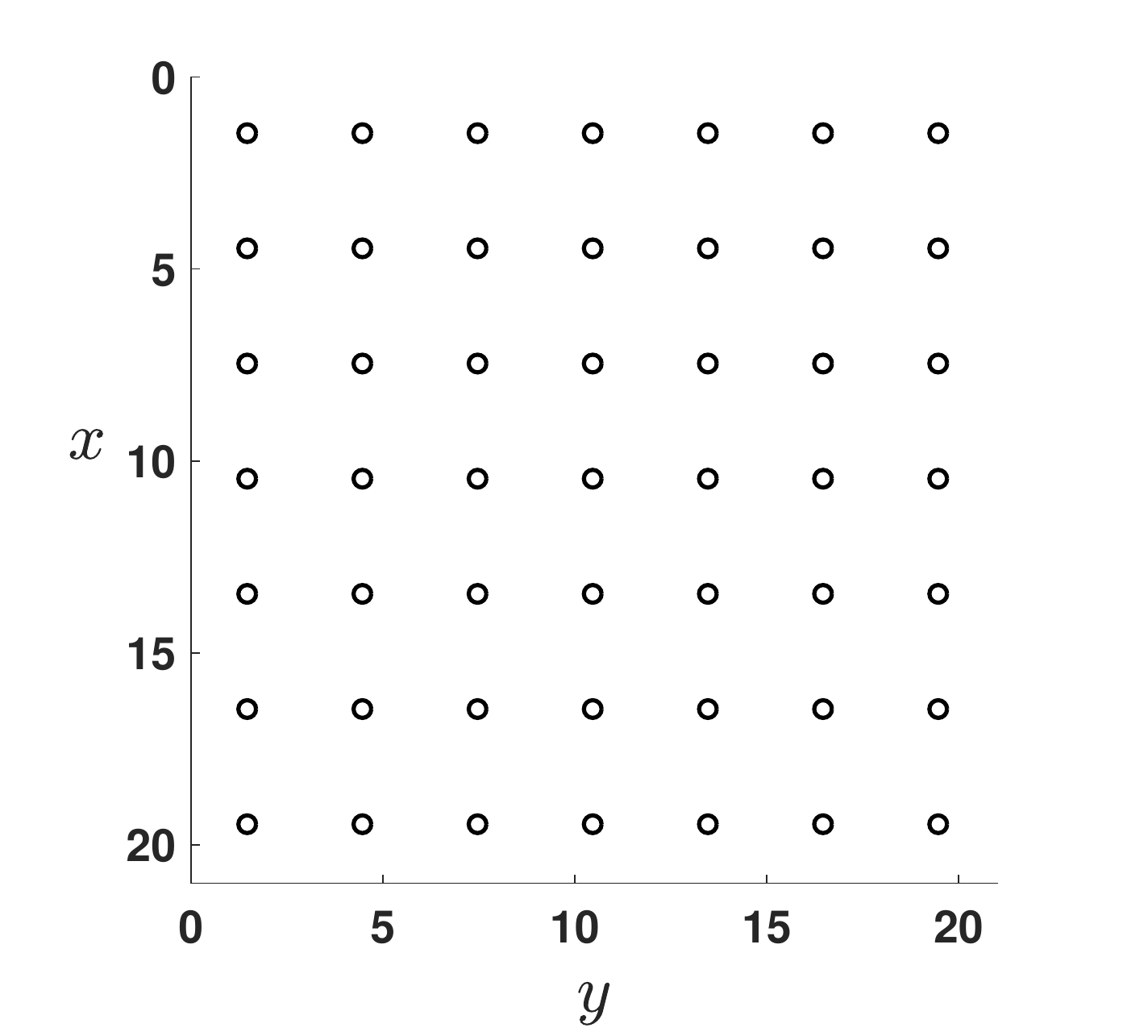}
\vspace{0.2cm}
\end{subfigure}
\begin{subfigure}{2.8in}
\caption{Perturbed equidistant.}
\includegraphics[width=2.8in]{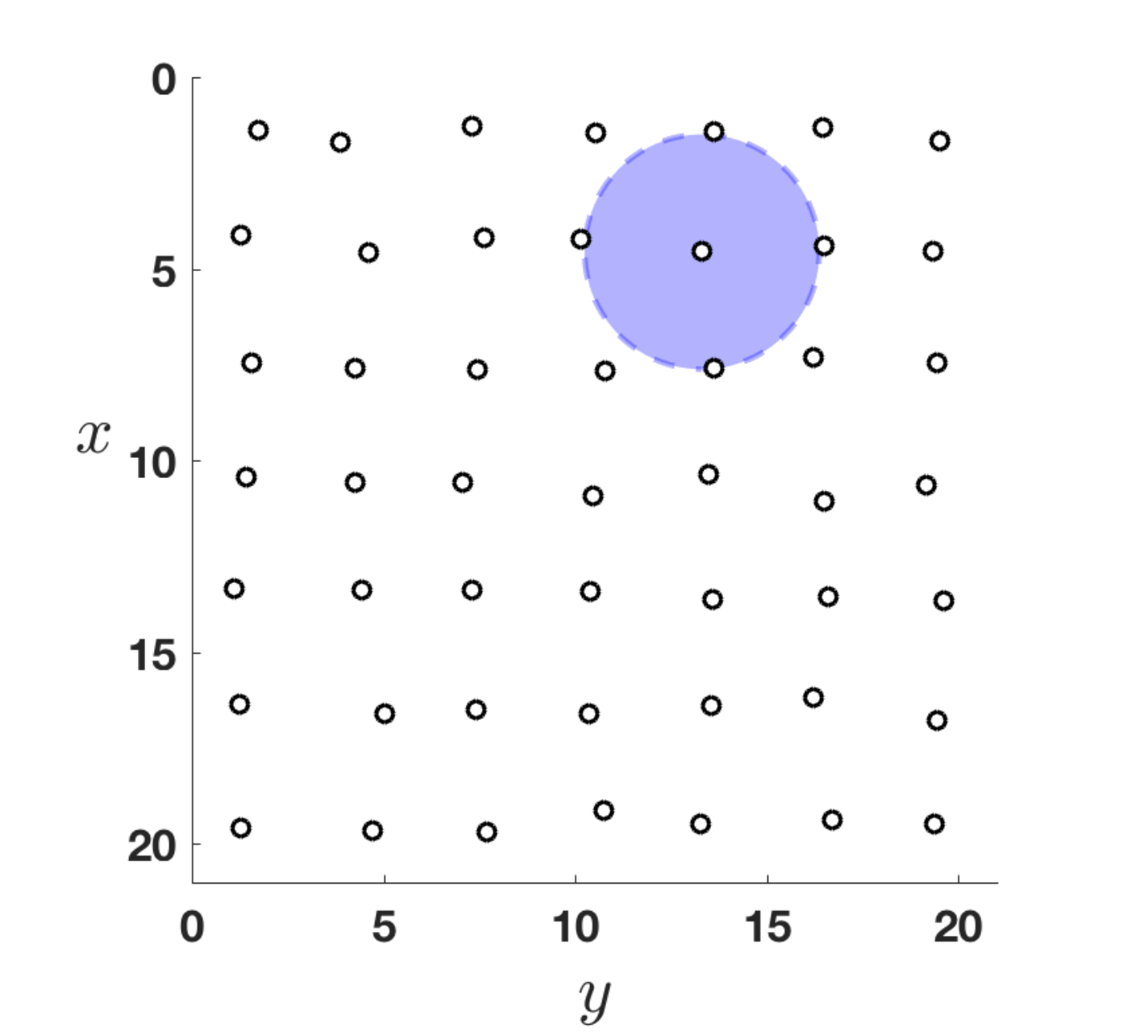}
\vspace{0.2cm}
\end{subfigure}
\begin{subfigure}{2.8in}
\caption{Random.} 
\includegraphics[width=2.8in]{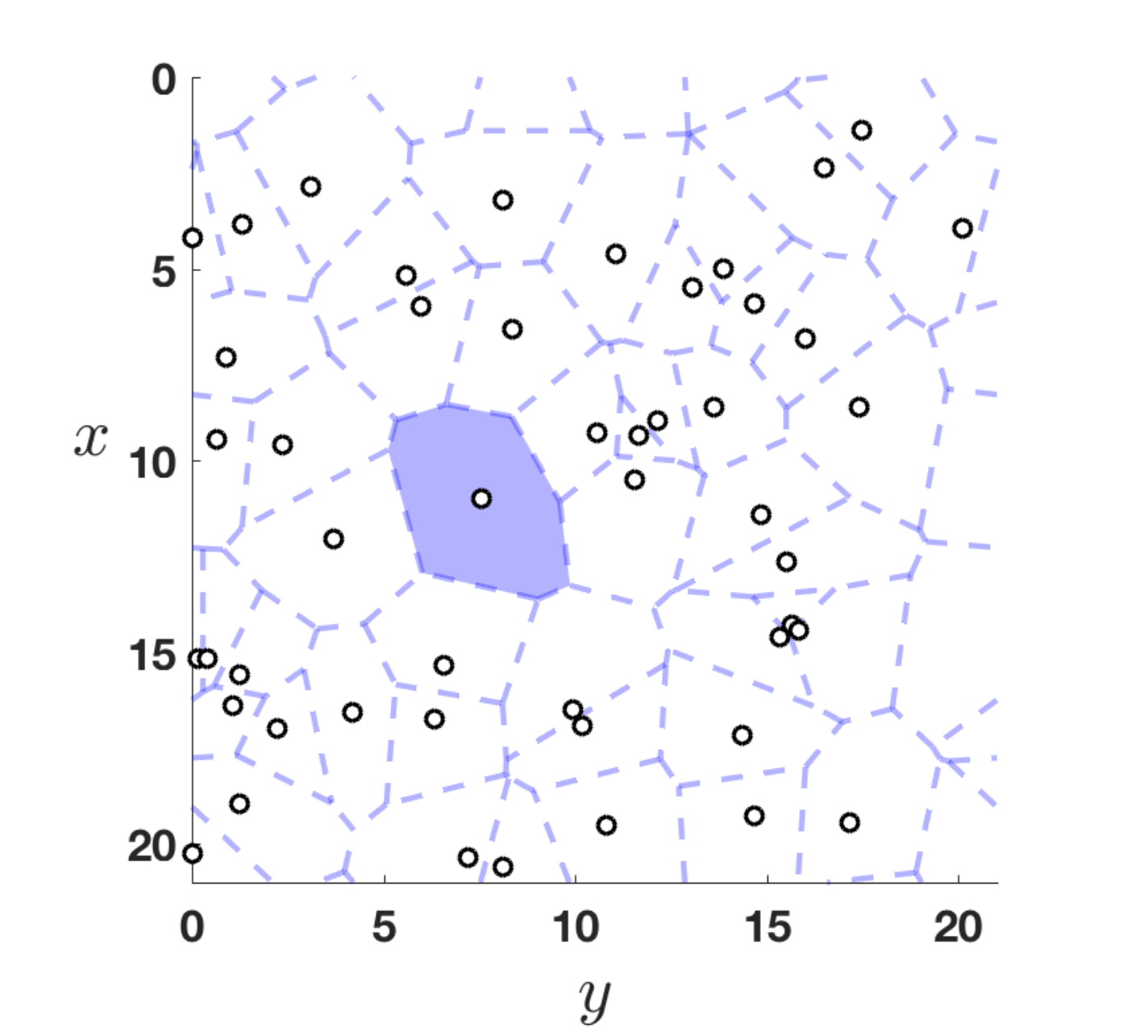}
\end{subfigure}
\begin{subfigure}{2.8in}
\caption{Quasirandom.}
\includegraphics[width=2.8in]{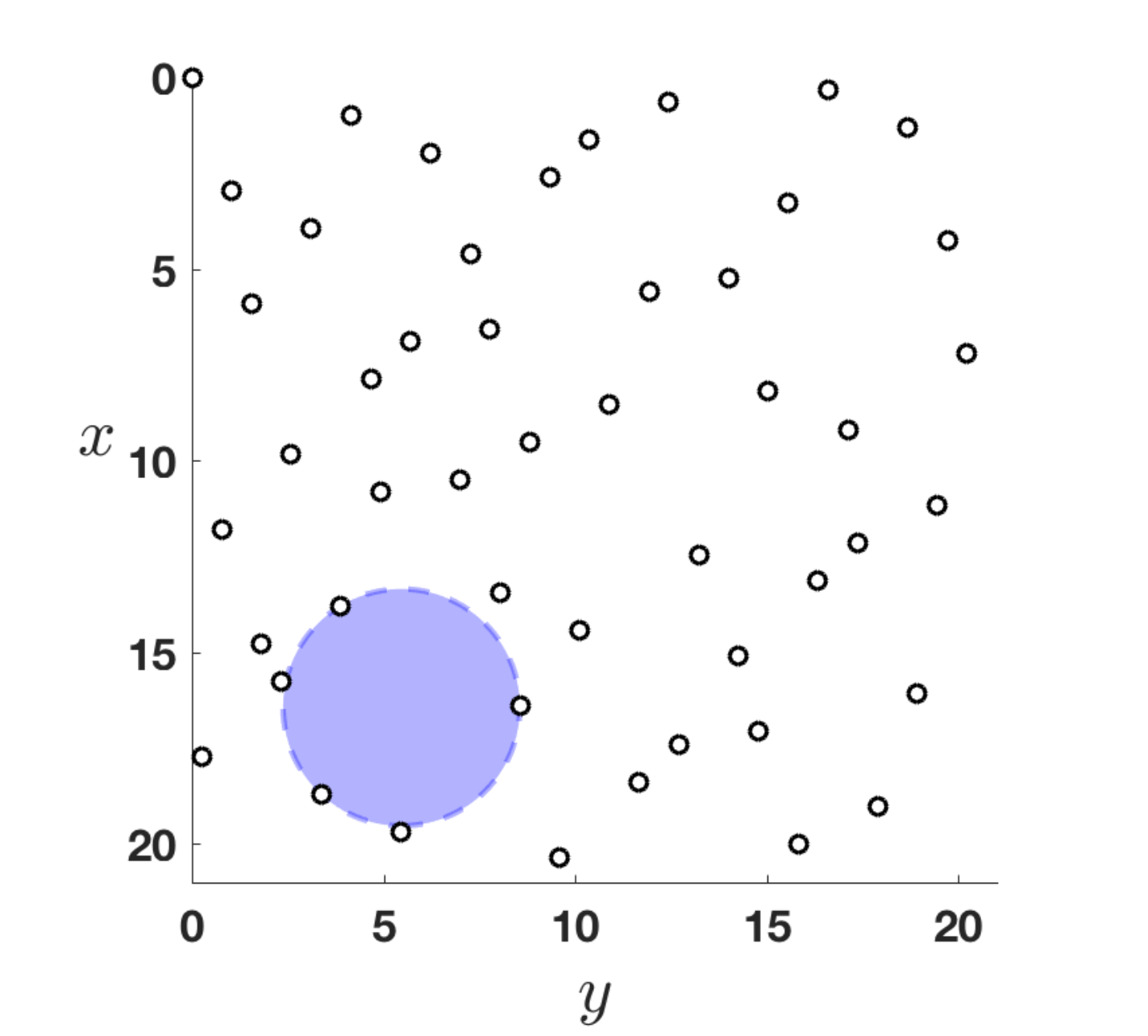}
\end{subfigure}
\caption{Examples of actuator grids used when $L_1 = L_2 = 21$ and $N_{\textrm{ctrl}} = 49$. The shaded regions shown in panels (b--d) correspond to calculations of the areas $A_1, A_2, A_3$ (defined in the text) which are used to measure the spread of the actuators.} \label{Examplegrids}
\end{figure}
In this paper, we consider a variety of actuator (and observer) arrangements. In 1D, \cite{gomes2016stabilizing} studied the optimal actuator placement for the 1D KSE \eqref{1dintroks}. For a given initial condition and finite time interval, they obtained the optimal (in terms of some cost functional) open-loop control with contributions from a finite set of optimally placed point actuators. The control locations were heavily dependent on the initial condition, desired state, and time interval over which controls were applied, thus the results do not imply that any particular arrangement of actuators performed well across a broad range of initial conditions. In this work, we do not seek results that are initial condition dependent or optimised for a cost functional, instead seeking actuator arrangements that give the best control performance in general. We utilise the following families of grids in our numerical simulations, examples of which are shown in Figure \ref{Examplegrids}:
\begin{enumerate}[(a)]
\item {\bf{Equidistant.}} These grids comprise of actuators which are equally spaced in the $x$- and $y$-directions with separations $d_1$ and $d_2$, respectively, forming a rectangular lattice (or a square lattice if $d_1 = d_2$). In order to comply with the periodicity of the domain, both $L_1/d_1$ and $L_2/d_1$ must be positive integers.
\item {\bf{Perturbed equidistant.}} An equidistant grid with actuators that have been randomly shifted by an amount smaller than the grid spacing (so that one actuator remains in each $d_1 \times d_2$ rectangular region). We sample the random shift from a normal distribution with zero mean.
\item {\bf{Random.}} The locations of the actuators are obtained by sampling coordinates from uniform distributions, i.e. $x_j \sim {\it{Unif}}(0,L_1)$, $y_j \sim {\it{Unif}}(0,L_2)$. Repeated locations are discarded.
\item  {\bf{Quasirandom.}} Also known as low-discrepancy sequences, quasirandom sequences are commonly used to sample space more evenly than uniform distributions. We use the 2,3-Halton sequence \citep{halton1960efficiency} to obtain quasirandom sequences of points in $[0,1]\times [0,1]$, which are appropriately rescaled to yield actuator locations in $Q$.
\end{enumerate}

Equally spaced (internal) actuators are commonly considered in studies of 1D control problems, however such arrangements are unsuccessful when the actuators are located at zeros of unstable eigenfunctions. Obviously there are actuator arrangements of interest not in the above classes which we do not study here, for example parallelogram lattice arrangements or hexagonal patterns. The success of a control strategy for a particular grid is equation and boundary condition dependent. The grids we consider appeared most natural for a study on rectangular periodic domains -- for a problem where periodicity is enforced on a hexagonal domain, a hexagonal actuator grid would be the most appropriate.

We measure three areas to quantify the spacing of the various grids -- examples of each are shown as shaded regions in Figure  \ref{Examplegrids}(b--d). The first, $A_1$, is defined as the area of the largest circle centred at an actuator which contains no other actuators, shown in panel (b) for the case of a perturbed equidistant grid. If all actuators are placed in groups of more than one, then $A_1$ will be small, with an infimum of zero found in the limit of each actuator approaching another. The supremum of $(L_1^2 + L_2^2)\pi/4$ corresponds to moving one actuator away from a cluster of all the other actuators. The value of $A_1$ is thus not entirely informative, but we find that its deviation from the value for the equidistant case, $A_1^E = \pi \cdot \min\{d_1,d_2\}^2$, is a more appropriate metric for our study. The quantity $A_2$ is defined via the Voronoi tessellation, shown in panel (c) with dashed lines, which separates the plane into the sets of points which are closest in Euclidean distance to each actuator. We define $A_2$ to be the area of the largest Voronoi cell, with the minimum value attained for the equidistant case, $A_2^E = d_1d_2$. Lastly, $A_3$, shown in Figure \ref{Examplegrids}(d), is defined as the area of the largest circle that can be inscribed in the plane without containing any actuators -- this is the solution of the well-known largest empty circle problem \citep{Toussaint1983}. This last area is minimised for equidistant actuator grids with $A_3^E = \pi (d_1^2 + d_2^2)/4$, and in contrast to $A_1$ and $A_2$, quantifies the size of the gaps between the actuators.

\section{Proportional control\label{SecPropCont}}

In this section, we consider proportional point-actuated controls. We assume that both observation and actuation occur at the same set of locations in the periodic domain $Q$. The strength and orientation of actuation at a point is based only on the observation of the local interface height at that instant in time. There is no communication from observers at other actuator locations, and no \textit{a priori} knowledge of the governing equation is utilised; this is the most basic level of feedback control.

To motivate our study, we allow actuation and observation at every location in space, and consider controls of the form $\zeta = -\alpha\eta$ for a strength $\alpha \geq 0$ to stabilise the trivial zero solution, i.e. $\overline{\eta} = 0$. The dispersion relation for the controlled KSE \eqref{controlled2dks} is then
\begin{equation}\label{lineargrowthrate}s(\bm{\tilde{k}}) = (1- \kappa)\tilde{k}_1^2 - \kappa \tilde{k}_2^2 - | \bm{\tilde{k}}|^4 - \alpha,\end{equation}
where $s$ is the linear growth rate and $\bm{\tilde{k}}$ is the scaled wavenumber vector for $\bm{k} \in \mathbb{Z}^2$ given by \eqref{ktildedefn}. Multiplying \eqref{controlled2dks} by $\eta$ and integrating by parts, we find that the energy of the solution (given by the $L^2$-norm) 
evolves according to the energy equation
\begin{equation}\frac{1}{2}\frac{\mathrm{d}}{\mathrm{d}t} \| \eta \|_{L^2}^2 = \sum_{\bm{k} \in \mathbb{Z}^2} s(\bm{\tilde{k}}) |\eta_{\bm{k}}|^2.\end{equation}
If $s < 0$ for all arguments, the $L^2$-norm of $\eta$ decays exponentially for all choices of $L_1$ and $L_2$; it follows from \eqref{lineargrowthrate} that this is achieved for $\alpha > \alpha_{\textrm{c}}(\kappa)$, where
\begin{equation}\label{alphacrit}\alpha_{\textrm{c}}(\kappa) = \begin{cases}
   \quad\;\; 0       & \quad \text{if } \kappa > 1,\\
   \displaystyle\frac{(1-\kappa)^2}{4} & \quad \text{if } \kappa \leq 1.\\
  \end{cases}
\end{equation}
With periodic boundary conditions, taking $\alpha > \alpha_{\textrm{c}}$ is sufficient, but a sharp bound may be obtained in terms of $L_1$ and $L_2$. For a non-trivial desired state, this control generalises to $\zeta = - \alpha ( \eta - \overline{\eta})$, and condition \eqref{alphacrit} becomes more complicated since the system must be linearised about a non-trivial state.

To explicitly obtain a critical actuation strength we used knowledge of the governing equation; if a PDE is well-posed in the very weak sense that its linear part is dissipative at small scales, then there exists an $\alpha_{\textrm{c}}$ such that for all $\alpha > \alpha_{\textrm{c}}$, the linear controlled system is stable. For an unknown system, this may be found experimentally by simply increasing the actuation strength $\alpha$. It is not necessarily true that this will yield full nonlinear stability of the given system, however, we expect that many interfacial problems with weak nonlinear interactions may respond well to this form of control.

It may not be viable to actuate at every spatial location, or even observe the entire interface. For the remainder of the section, we consider point-actuated controls \eqref{PAcontrolform1} with
\begin{equation}\label{phiform1}\phi^{j}(t;\eta, \overline{\eta}) = - \alpha [\eta(\bm{x}_{j},t) - \overline{\eta}(\bm{x}_{j},t)],\end{equation}
where $\alpha \geq 0$ is the actuation strength -- this is often referred to as ``pinning control" in the literature \citep{grigoriev1997pinning}. We note that the solution average is not necessarily preserved by this choice of forcing, but successful controls ensure that solutions do not deviate much from zero mean.

\subsection{Convergence study}

For our convergence study, we take $L_1 = L_2 = 21$ and $\kappa = 0.25$ (corresponding to an overlying film). We use $49$ actuators/observers in a quasirandom grid which are located according to the 2,3-Halton sequence -- the arrangement is shown in Figure \ref{Examplegrids}(d). These parameters will be employed in a later subsection for a comparison of grids. Here, we fix the initial condition to be
\begin{align}\label{initicond1}\eta_0(\bm{x}) = \frac{1}{10} & \left[ \cos\left(\frac{2\pi x}{L_1} \right) + \cos\left(\frac{2\pi x}{L_1} + \frac{2\pi y}{L_2}  \right) \right. \nonumber \\
&\qquad\qquad\quad\;\; \left. + \sin\left(\frac{4\pi x}{L_1} + \frac{2\pi y}{L_2}  \right) +  \sin\left( \frac{2\pi y}{L_2}  \right) + \sin\left( \frac{4\pi y}{L_2}  \right) \right],\end{align}
and let the system evolve until $t=1$ without controls. After this time, we actuate proportionally with strength $\alpha = 150$ to drive the solution to the zero state, $\overline{\eta} = 0$, until $t = 2$. In the following table we give the value of the $L^2$-norm at time $t=2$, $\mathcal{C}_1(2)$, for a range of spatial discretisations $M$, $N$, and time discretisation $\Delta t$:
\begin{table}[H]
\centering
\begin{tabular}{|c|>{\hspace{-10pt}}c<{\hspace{-10pt}}|>{\hspace{-10pt}}c<{\hspace{-10pt}}|>{\hspace{-10pt}}c<{\hspace{-10pt}}|}\hline
\backslashbox{$\Delta t$}{$M = N$}& $128$ & $256$ & $512$\\
\hline
$2 \times 10^{-4}$ & --- & --- & --- \\ \hline
$1 \times 10^{-4}$ & $0.045441843348824$ & $0.045456708994765$ & $0.045460357974812$ \\ \hline
$5 \times 10^{-5}$ & $0.045442820551101$ & $0.045457704277985$ & $0.045461333339984$ \\ \hline
$2.5 \times 10^{-5}$ & $0.045443306597597$ & $0.045458189680699$ & $0.045461847168330$ \\ \hline
$1 \times 10^{-5}$ & $0.045443600470486$ & $0.045458433050053$ & $0.045462134162334$\\ \hline
\end{tabular}
\caption{Values of $\mathcal{C}_1(2)$ for a range of spatial and temporal discretisations.}\label{tab:convergence}
\end{table}
Recall that we are employing a fourth-order BDF scheme which is not unconditionally stable, unlike the first- and second-order BDF schemes. In a convergence study of the unforced equation, \cite{akrivislinearly} observed, for a particular set of parameters, that once the higher-order BDF schemes (third- to sixth-order) are stable, the solution converges to machine accuracy (even for the worst case in Table \ref{tab:convergence} we found that $\mathcal{C}_1(1) = 3.265272$ was accurate to $8$ significant figures). Analogously, with the addition of point actuators, we find that once stability of the fourth-order scheme has been achieved (the scheme is not stable for $\Delta t = 2\times 10^{-4}$), then $\mathcal{C}_1(2)$ is accurate to $3$ decimal places ($2$ significant figures). For the unforced equation, numerical solutions showed that the spectrum of solutions decays exponentially, indicating analyticity \citep{tkp2018}. It was also observed that the spectrum decays faster in the transverse $\tilde{k}_2$ modes than the streamwise $\tilde{k}_1$ modes; this is expected given the asymmetry of the linear part of \eqref{controlled2dks}. In numerical simulations on a square periodic domain, the Fourier mode truncation $N$ is not required to be as large as $M$ for good accuracy. This is no longer true with the introduction of Dirac delta functions -- the decay of the Fourier spectrum becomes more symmetric in wavenumber space, placing similar restrictions on the constants $M$ and $N$ to obtain good accuracy. Thus, in Table \ref{tab:convergence}, we only consider $M = N$. 
\begin{figure}
\centering
\includegraphics[width=3.5in]{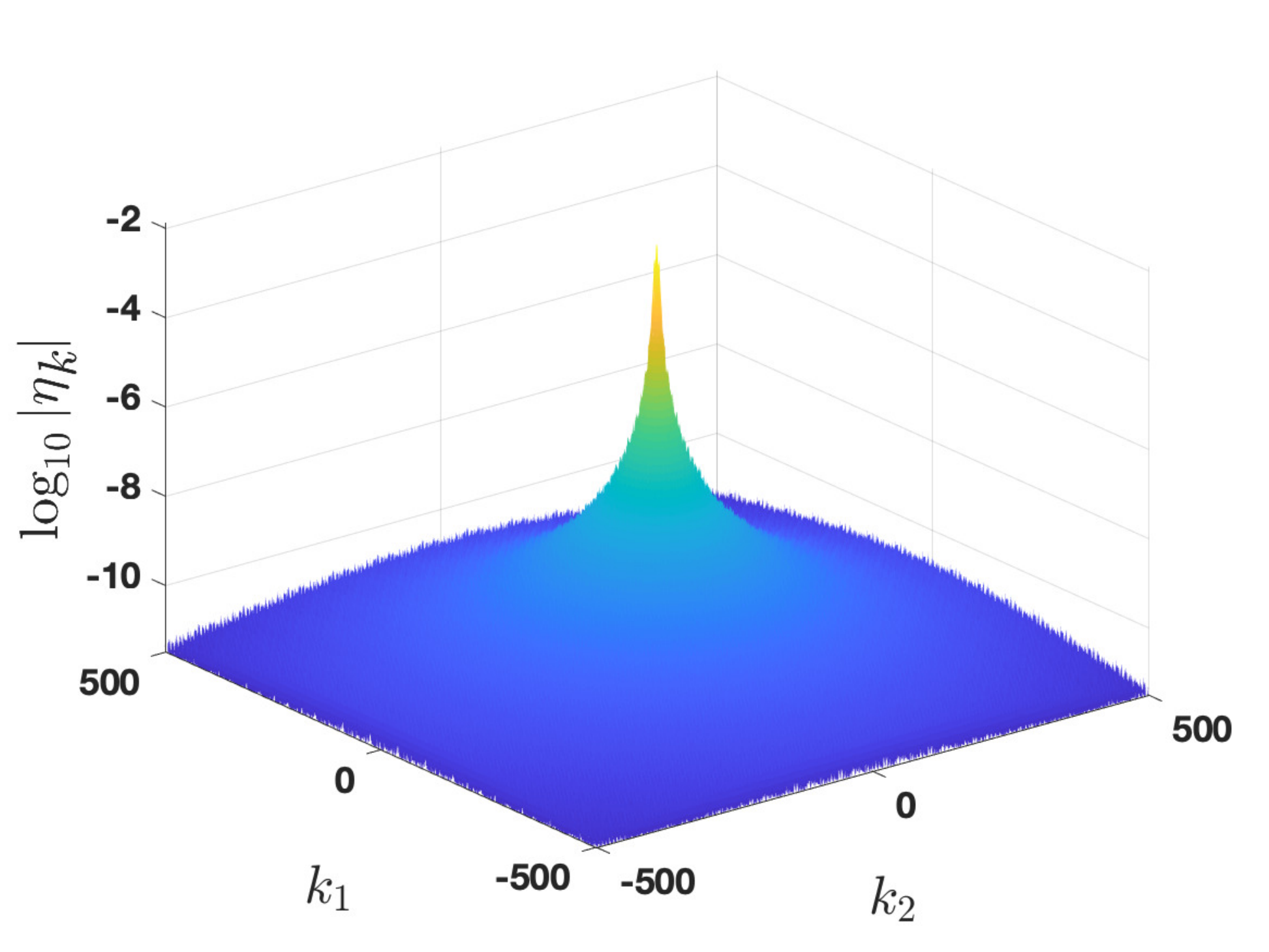}
\caption{Spectrum of the controlled solution at time $t=2$.}\label{Fufinalconvergencetest}
\end{figure}
The spectrum at $t = 2$ from the most well-resolved simulation is shown in Figure \ref{Fufinalconvergencetest} (plotted against the unscaled wavenumbers $|k_1|,|k_2| \leq 500$). Analyticity (exponential decay of the spectrum) is lost due to the Dirac delta forcing, and we find numerically that $|\eta_{\bm{k}}| \sim |\bm{\tilde{k}}|^{-4}$. Thus, $\Delta^2 \eta$ has an approximately constant spectrum, $|\bm{\tilde{k}} |^4|\eta_{\bm{k}}| \sim O(1)$,  which balances the constant spectrum of the Dirac delta forcing.

\subsection{Controlling unbounded exponential growth\label{subsectionsuppressgrohanging}}

% This figure can be made using /Users/rubentomlin/Dropbox/1st Year PhD/2DControlruns/Propcontrolhanging/six

\begin{figure}
\centering
\begin{subfigure}{2.8in}
\caption{Costs for $N_{\textrm{ctrl}}=100$.} 
\includegraphics[width=2.8in]{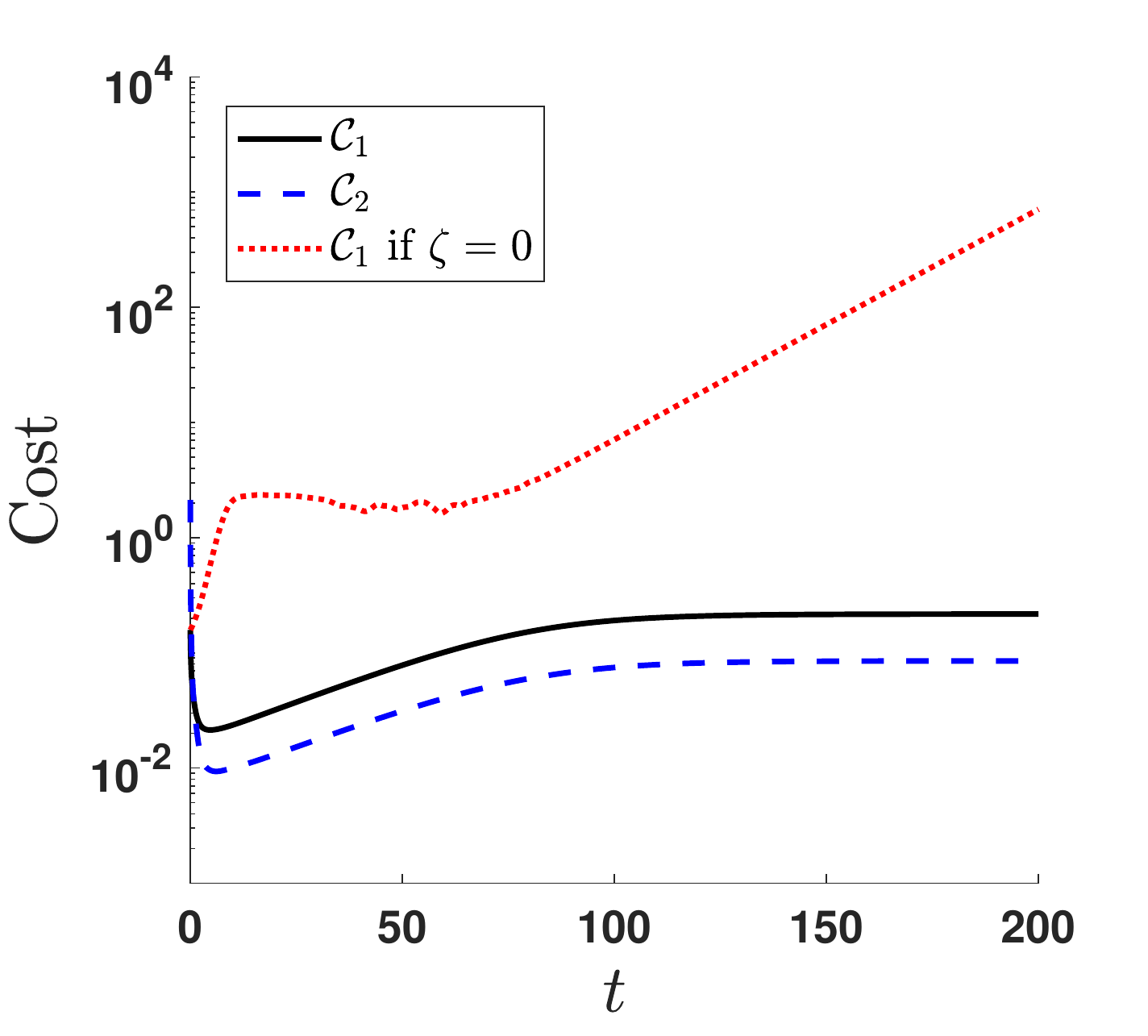}
\vspace{0.2cm}
\end{subfigure}
\begin{subfigure}{2.8in}
\caption{Solution at $t=200$ for $N_{\textrm{ctrl}}=100$.} \label{hangingcontourplat}
\includegraphics[width=2.8in]{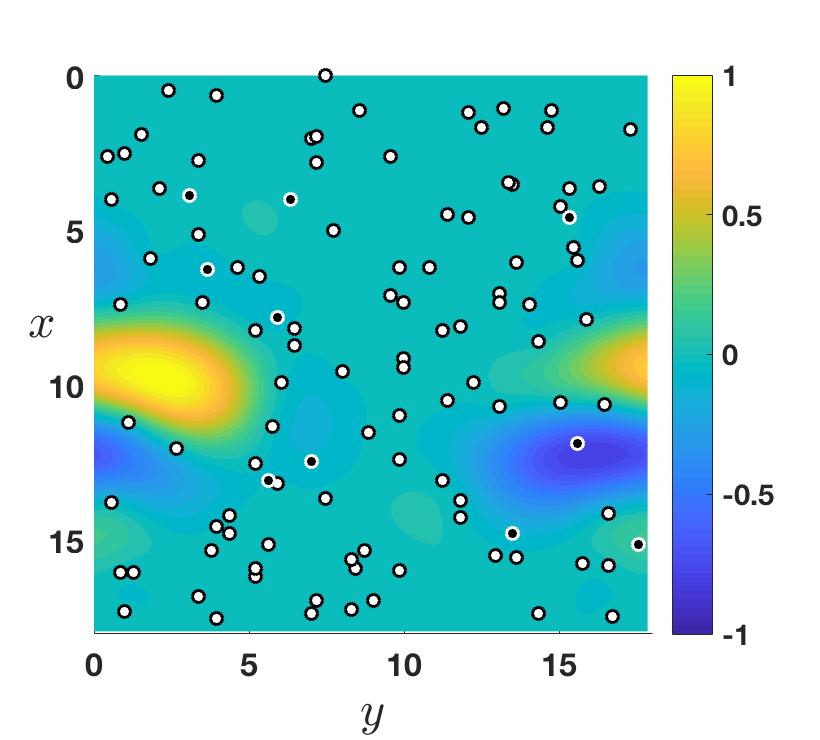}
\vspace{0.2cm}
\end{subfigure}
\begin{subfigure}{2.8in}
\caption{Costs for $N_{\textrm{ctrl}}=110$.} 
\includegraphics[width=2.8in]{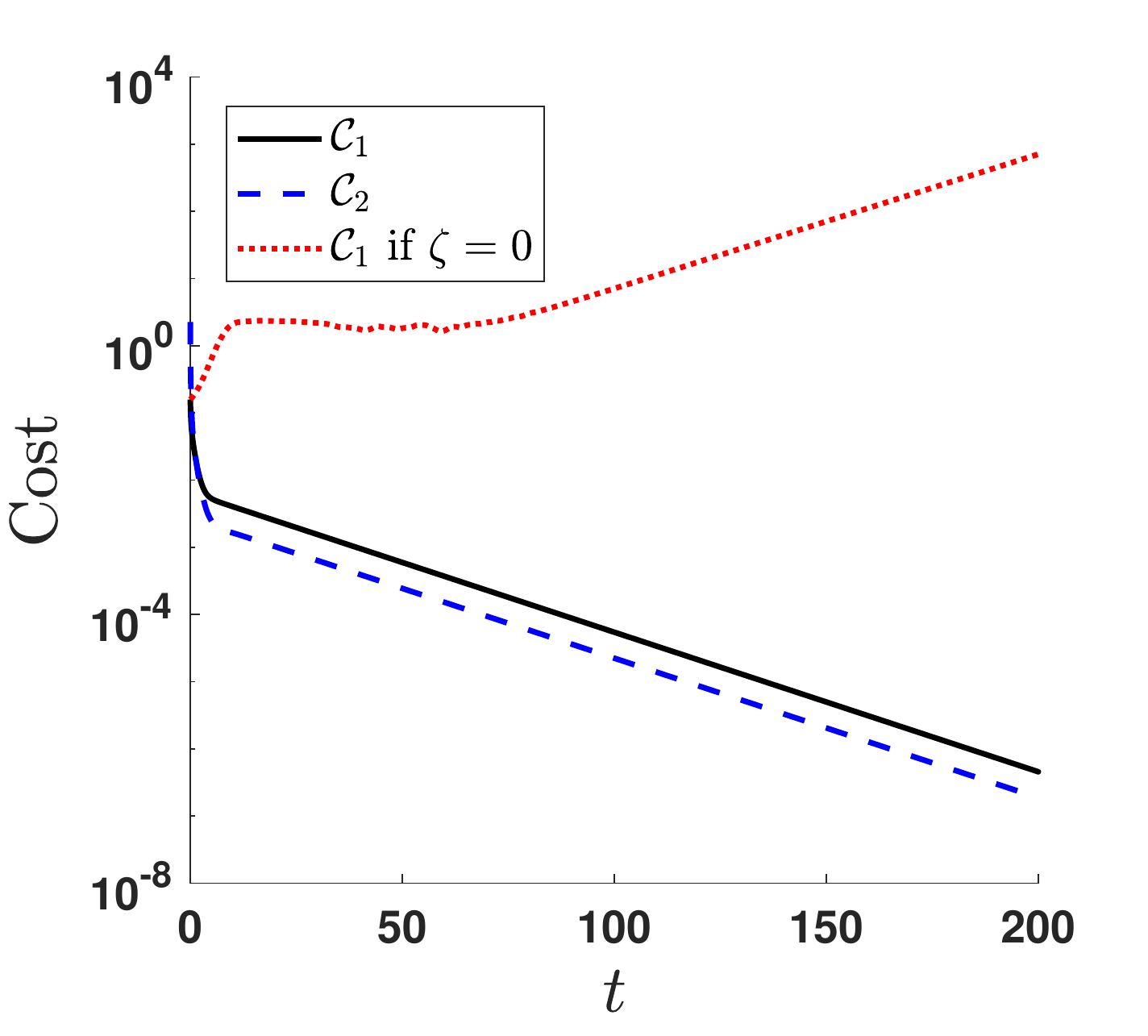}
\end{subfigure}
\begin{subfigure}{2.8in}
\caption{Solution at $t=200$ for $N_{\textrm{ctrl}}=110$.} \label{hangingcontourexp}
\includegraphics[width=2.8in]{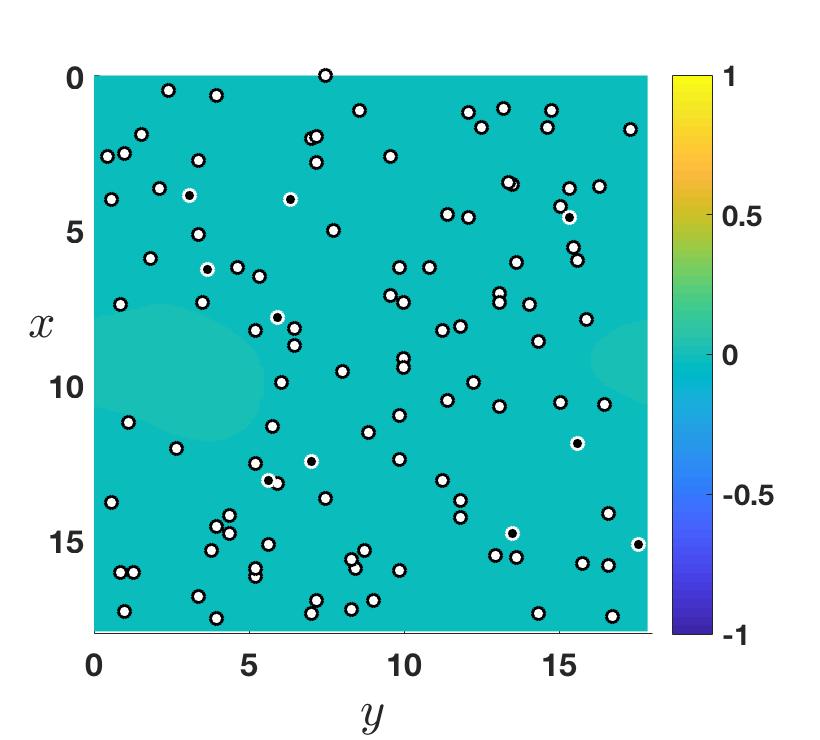}
\end{subfigure}
\caption{Costs and solution contours for $\alpha = 55$ and $N_{\textrm{ctrl}}=100, \; 110$. Panels (a) and (c) display the evolution of the costs $\mathcal{C}_1$ and $\mathcal{C}_2$ for the cases of $N_{\textrm{ctrl}}=100$ and $N_{\textrm{ctrl}}=110$, respectively. Panels (b) and (d) display contours of the respective solution profiles at $t=200$. The actuator locations are superimposed on the solution contours in panels (b) and (d) -- the white circles with black edges denote the point actuators used in both cases, and the black circles with white edges denote the remaining $10$ which are ``switched on" in the $N_{\textrm{ctrl}}=110$ case.} \label{hangingpropcostsC1C2}
\end{figure}

\begin{figure}
\centering
\begin{subfigure}{2.8in}
\caption{Regimes in the $\alpha$--$N_{\textrm{ctrl}}$ plane.} 
\includegraphics[width=2.8in]{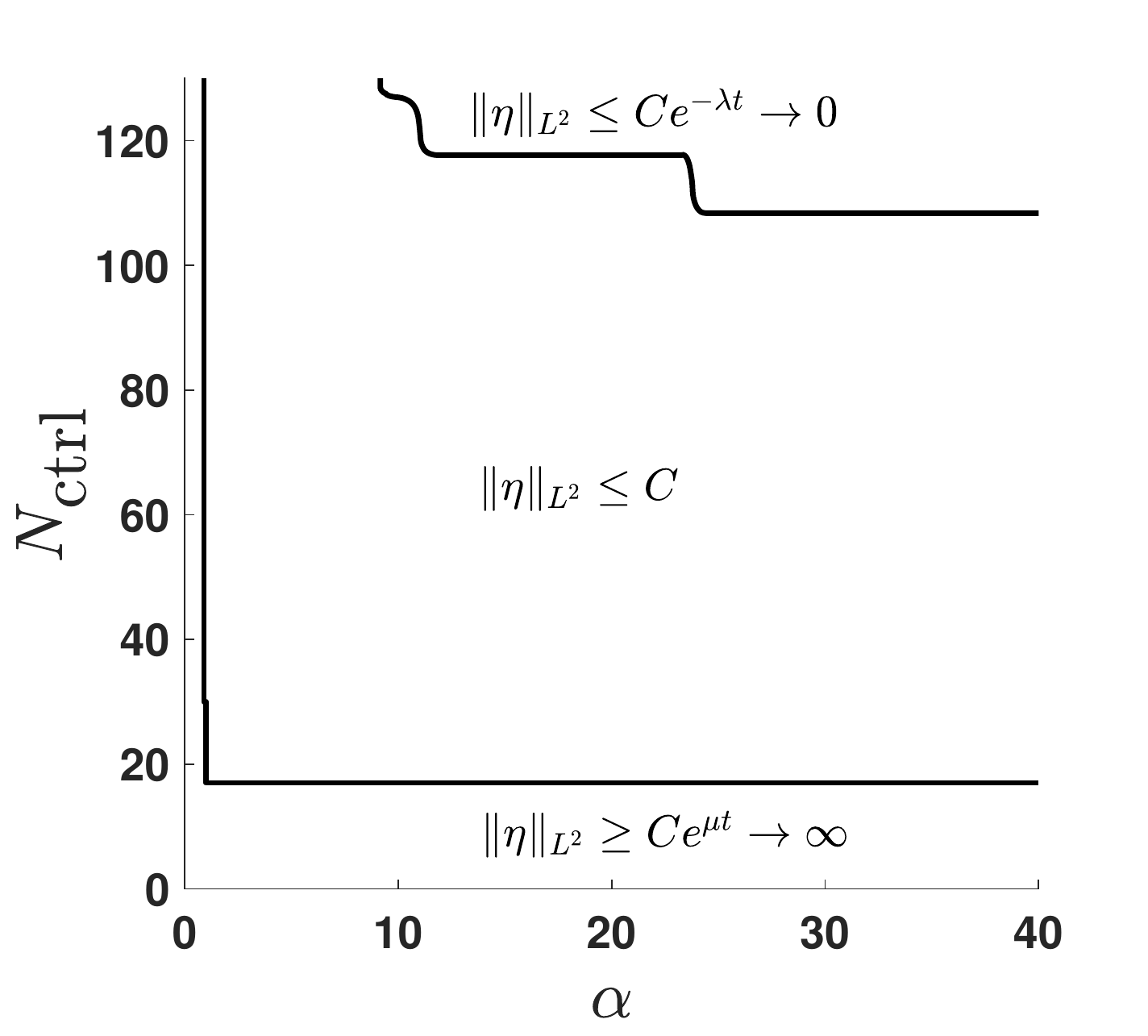}
\end{subfigure}
\begin{subfigure}{2.8in}
\caption{Decay rate $\lambda$ for $N_{\textrm{ctrl}}=110$.} 
\includegraphics[width=2.8in]{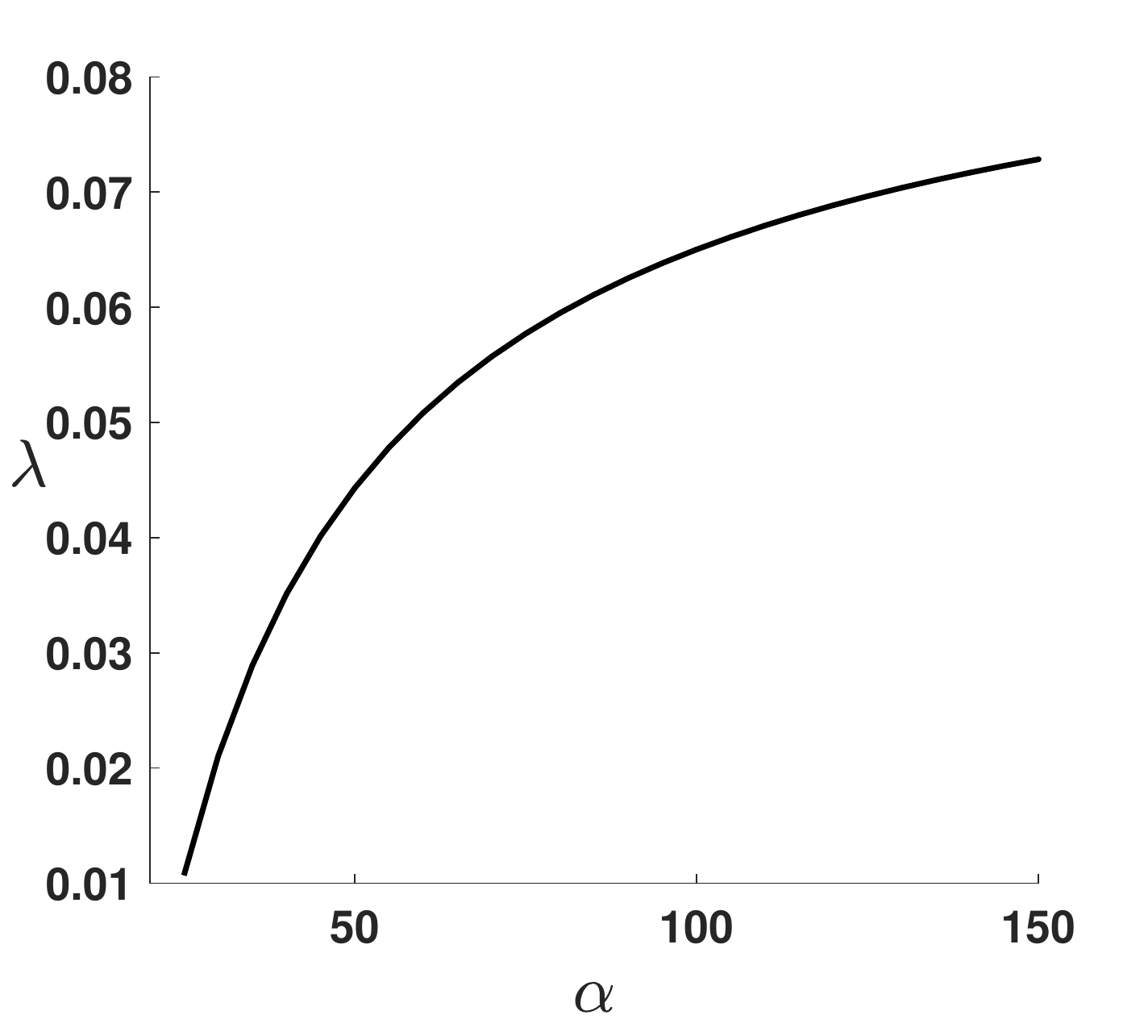}
\end{subfigure}
\caption{Success of controls for hanging films as $\alpha$ and $N_{\textrm{ctrl}}$ vary. Panel (a) displays three regimes dependent on the values of the control parameters. The regions shown are based on over 500 numerical simulations. Panel (b) shows the exponential decay rate $\lambda$ of the costs $\mathcal{C}_1$ and $\mathcal{C}_2$ as a function of $\alpha$ for $N_{\textrm{ctrl}}=110$.}\label{stabofhangingfilmfig}
\end{figure}

In this subsection, we investigate the efficacy of proportional point-actuated controls in suppressing the unbounded growth of solutions to the 2D KSE \eqref{controlled2dks} with $\kappa < 0$; we also set $\overline{\eta} = 0$. This choice of $\kappa$ physically corresponds to a hanging fluid film, and the linear instabilities in the transverse modes are the classical Rayleigh--Taylor instabilities of fluid dynamics. Additionally, there are the usual instabilities in the streamwise and mixed modes present for $0 \leq \kappa < 1$. Optimal controls involving the transverse modes alone were applied to \eqref{controlled2dks} in \citep{tomlin2019optimal}, revealing windows of chaotic and travelling wave attractors for the streamwise and mixed modes. The success of the controls will be measured by two objectives; the first being the successful suppression of transverse growth resulting in bounded solutions, and the second being the stronger property of exponential stability of the flat film solution.

We fix $\kappa = -0.5$ and consider a square domain with $L_1 = L_2 = 18$, for which $30$ Fourier modes are linearly unstable in total, including the transverse $(0,\pm 1)$- and $(0,\pm 2)$-modes. For our numerical simulations, we take initial condition \eqref{initicond1}, having contributions from all unstable transverse modes. Controls are applied from the initial time, rather than letting the transverse instabilities develop further. For this subsection, we consider random grids which are constructed recursively as follows. Actuator locations are obtained by sampling coordinates from a uniform distribution, and are ordered in a list so that taking $N_{\textrm{ctrl}} = p$ corresponds to ``switching on'' the first $p$ control actuators/observers in the list. This arrangement and ordering is fixed in this subsection, across all simulations for different values of $\alpha$. Fixing the initial condition and having consistency in how the actuators are arranged and the order in which they are ``switched on'' allows us to draw robust conclusions. The success or failure of the controls is analysed in the $\alpha$--$N_{\textrm{ctrl}}$ plane.

Figure \ref{hangingpropcostsC1C2} shows the costs $\mathcal{C}_1$ and $\mathcal{C}_2$, and solution contours for two choices of $N_{\textrm{ctrl}}$ with $\alpha = 55$. Panels (a,b) correspond to the case $N_{\textrm{ctrl}} = 100$, where the costs plateau at a finite value and the solution contours show the presence of cellular humps in regions where no actuators are active. For the $N_{\textrm{ctrl}}=110$ case shown in panels (c,d), the zero solution is exponentially stabilised. The actuator locations for both cases are superimposed on the solution contours in panels (b) and (d), the black circles with white edges represent the actuators which are ``switched on" for $N_{\textrm{ctrl}}=110$ only. One of the actuators which is ``switched on" for the $N_{\textrm{ctrl}}=110$ case is located at a cellular hump which forms for the $N_{\textrm{ctrl}}=100$ simulation, indicating how crucial the location of point actuators may be, since cellular humps form in sufficiently large areas of $Q$ where no actuators are located. As discussed previously, given an initial condition, it is possible to optimise the actuator locations (see \cite{gomes2016stabilizing} for the 1D case). However, we are concerned with finding actuator spacings/arrangements which provide the best stabilisation of the zero solution for any initial condition -- this is investigated in the next subsection.

In order to determine how the control success depends on the parameters $\alpha$ and $N_{\textrm{ctrl}}$, numerical experiments were carried out with $N_{\textrm{ctrl}}$ ranging from $0$ to $130$, and $\alpha \in [0,150]$. Figure \ref{stabofhangingfilmfig}(a) provides the numerical results over a section of the $\alpha$--$N_{\textrm{ctrl}}$ plane, where the parameter space is seen to be divided into three regions depending on the behaviour of the costs $\mathcal{C}_1$ and $\mathcal{C}_2$. Unsurprisingly, for sufficiently small values of $\alpha$ and $N_{\textrm{ctrl}}$, the growth of the transverse modes is not saturated. For $\alpha \geq 1.5 $ and $N_{\textrm{ctrl}} \geq 18$, approximately, bounded solutions emerge as seen in Figure \ref{hangingpropcostsC1C2}(b). Exponential stabilisation is obtained for much larger values of the control parameters; we find that for $\alpha = 24$, the solution is bounded and non-zero for $N_{\textrm{ctrl}} = 108$, but find exponential stabilisation of the zero solution for $N_{\textrm{ctrl}}=109$. The boundaries of the region on the right of panel (a) continue at a constant value of $N_{\textrm{ctrl}}$ for the strengths $\alpha \in [40,150]$ not shown in the figure. The seemingly sharp critical values in both $\alpha$ and $N_{\textrm{ctrl}}$ are quite surprising; it would be more expected that the same level of control could be achieved for smaller $\alpha$ and increased $N_{\textrm{ctrl}}$ (or vice versa). In Figure \ref{stabofhangingfilmfig}(b), the exponential decay rate $\lambda$ of the costs $\mathcal{C}_1$ and $\mathcal{C}_2$ (the decay of these two quantities is approximately the same) is plotted against $\alpha$ for the grid with $N_{\textrm{ctrl}}=110$. It is evident that $\lambda$ is a monotonically increasing function of $\alpha$, but that there is a maximal exponential decay rate associated with each actuator grid. For Figure \ref{stabofhangingfilmfig}(b), the maximal value of $\lambda = 0.09$ is predicted by fitting in the limit of large strength $\alpha$.

\subsection{Comparison of actuator arrangements.\label{CompofGridssubsec}} 

\begin{figure}[h!]
\centering
\begin{subfigure}{2.8in}
\caption{Decay rate $\lambda$ against area $|A_1 - A_1^E|$}
\includegraphics[width=2.8in]{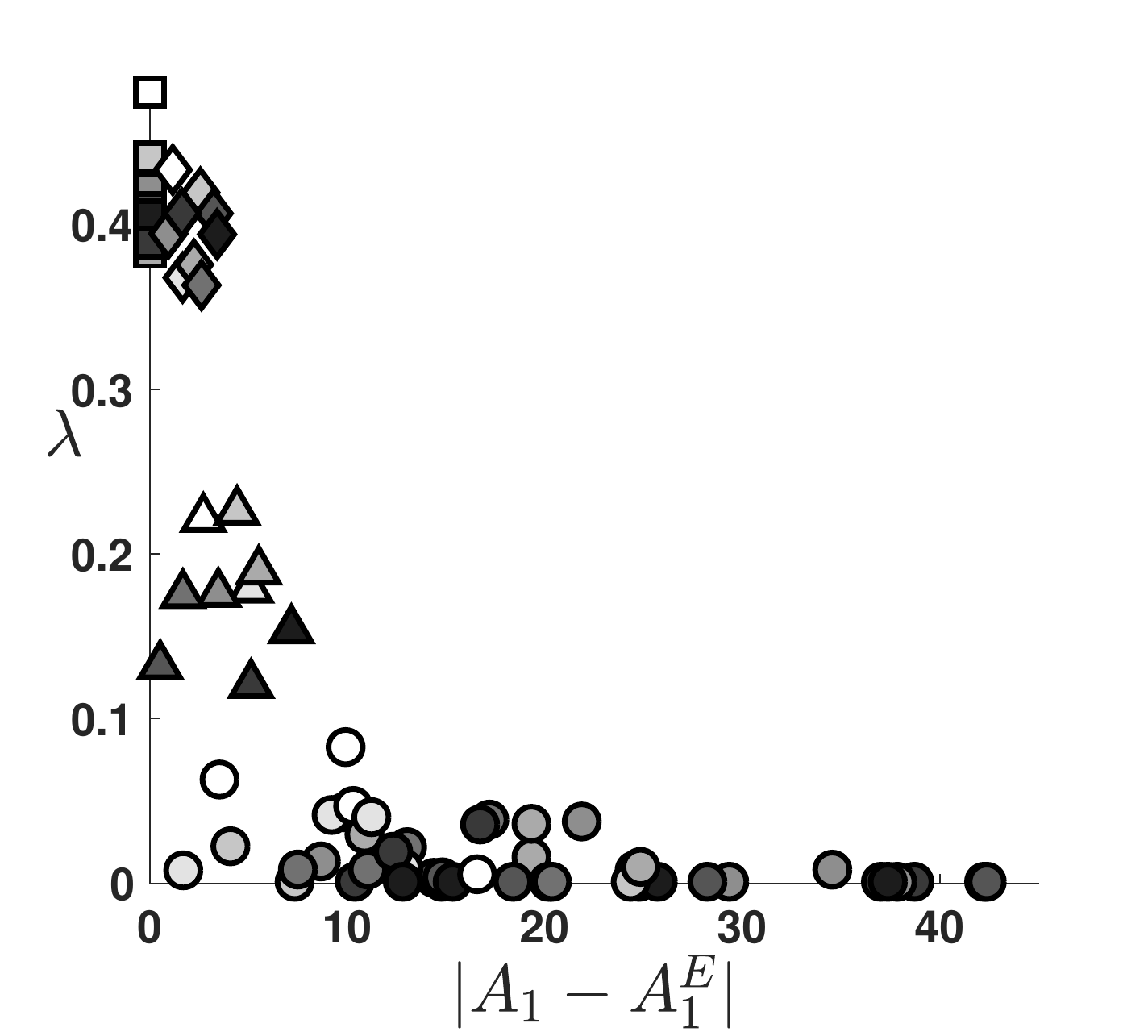}
\vspace{0.2cm}
\end{subfigure}
\begin{subfigure}{2.8in}
\caption{Decay rate $\lambda$ against area $A_2$}
\includegraphics[width=2.8in]{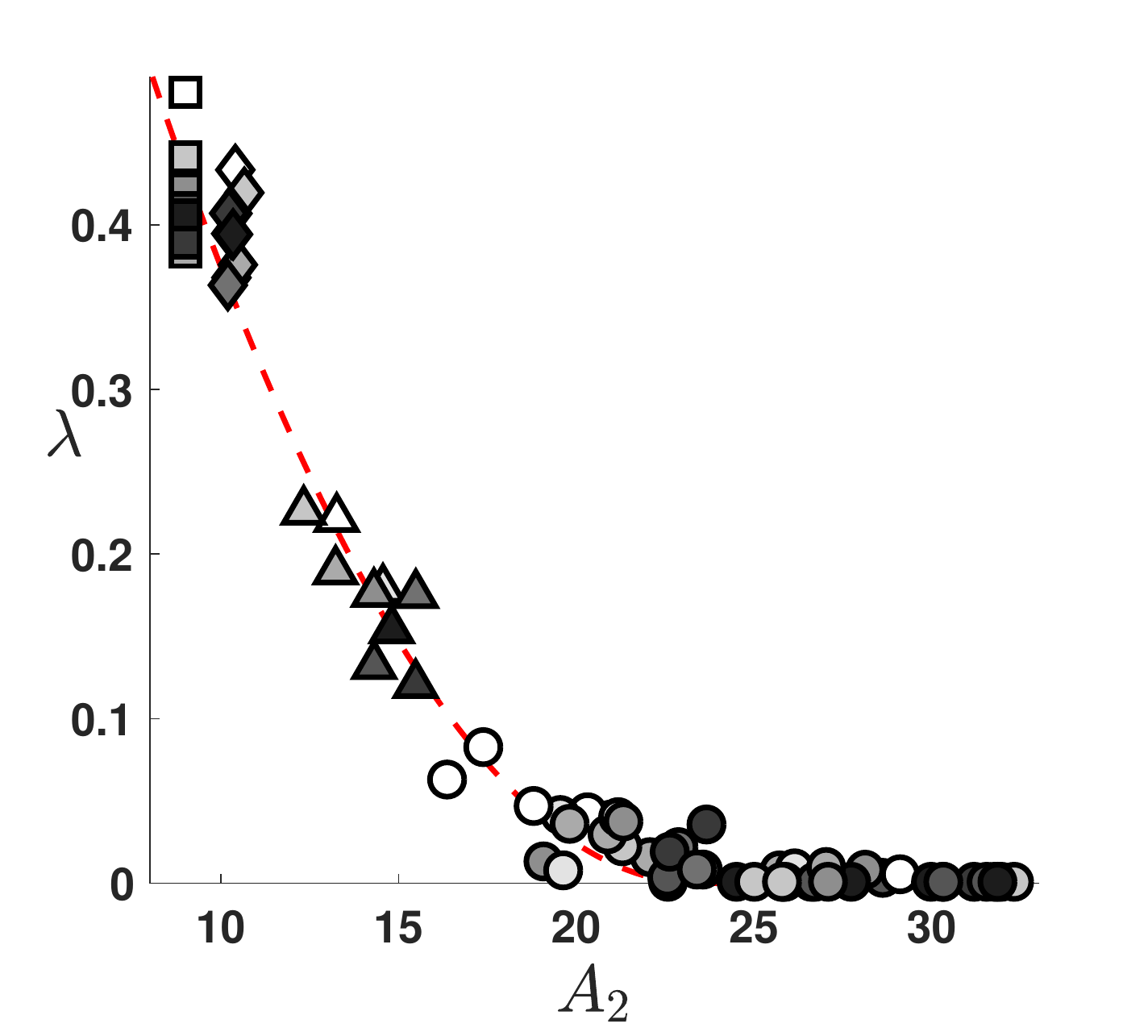}
\vspace{0.2cm}
\end{subfigure}
\begin{subfigure}{2.8in}
\caption{Decay rate $\lambda$ against area $A_3$}
\includegraphics[width=2.8in]{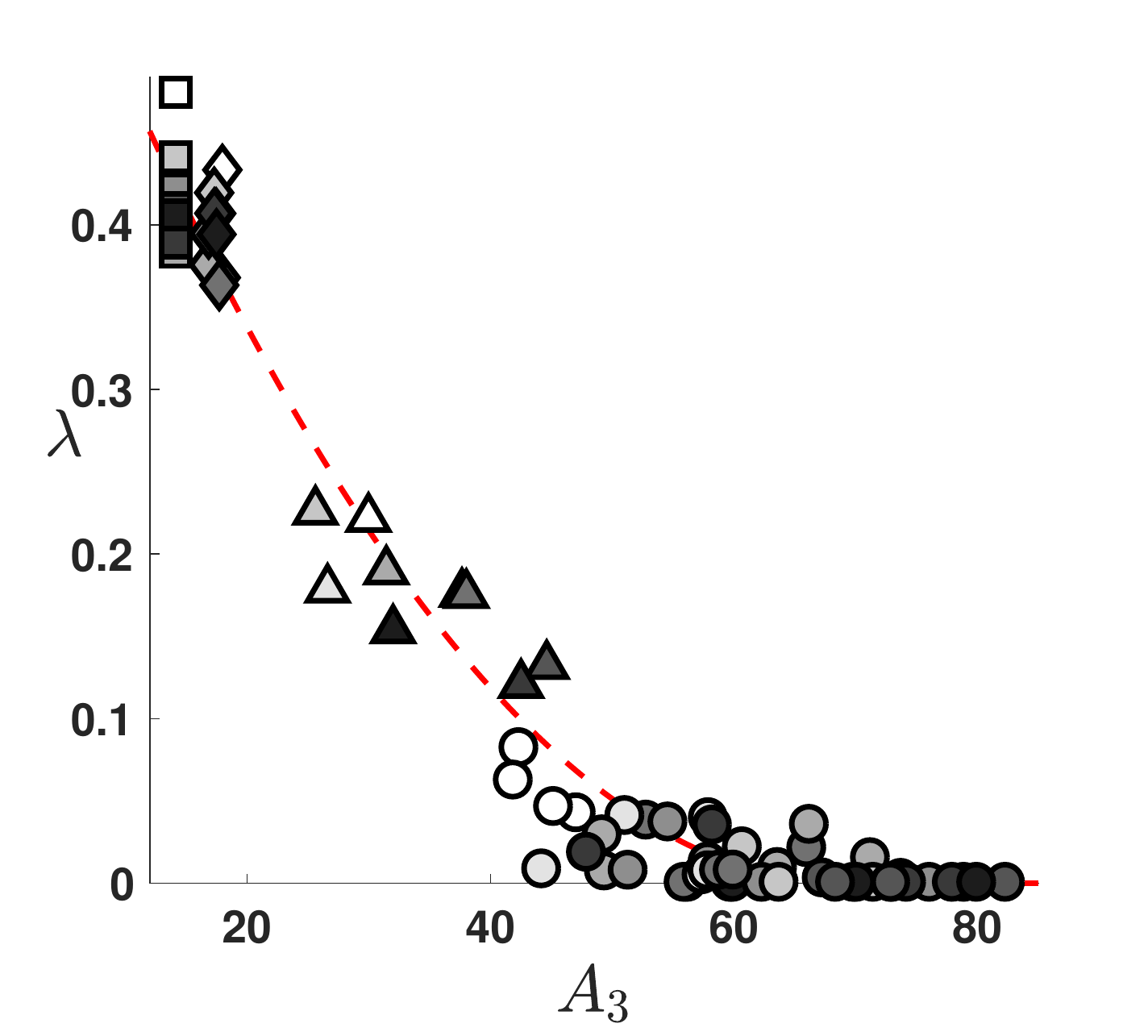}
\end{subfigure}
\caption{Performance of various actuator arrangements. The area predictors $|A_1 - A_1^E|$, $A_2$ and $A_3$ are plotted against the decay rate $\lambda$ of the costs $\mathcal{C}_1$ and $\mathcal{C}_2$. In cases when exponential decay is not observed, the value of $\lambda = 0$ is assigned. Markers: {\small $\square$} -- Equidistant; {\LARGE$\diamond$} -- Perturbed equidistant;  $\bigtriangleup$ -- Quasirandom; {\LARGE$\circ$} -- Random. The markers are shaded according to the values of $L$ in Table \ref{tab:numerics}, with the lightest corresponding to $L=21$ and the darkest to $L=45$. For $A_2$ and $A_3$, the data for $L=39,42,45$ is used to construct threshold models (details given in the text), given by the dashed lines in panels (b,c).} \label{propcontrolsupcritscatter}
\end{figure}

In this subsection, we test different actuator arrangements with the aim of understanding the relationship between actuator spacing and control performance. We vary the domain size $|Q|$, while keeping the ratio $|Q|/N_{\textrm{ctrl}}$ constant, i.e. we have a fixed number of actuators per unit area, and set $\alpha = 150$ (this is chosen to be large so that the control strength is not a limiting factor -- see the previous subsection). A finite energy density is observed for solutions of the uncontrolled system with $\kappa = 0$ on large domains \citep{tkp2018}, thus it is expected that the number of actuators required per unit area for successful control should be constant in the limit of large periodicities, $L_1, L_2 \rightarrow \infty$.

We set $\kappa = 0.25$, restrict to square domains with $ L_1 = L_2 = L$, and take $N_{\textrm{ctrl}}$ to be a square number so that we may compare the results of random and quasirandom arrangements with results for equidistant and perturbed equidistant actuator arrangements (with $d_1 = d_2$). The domain dimensions and numbers of actuators are chosen so that $|Q|/N_{\textrm{ctrl}} = 9$, i.e. one control actuator per 9 units of domain area, and are summarised in Table~\ref{tab:numerics}:
\begin{table}[H]
\centering
\begin{tabular}{| c | >{\hspace{-10pt}}c<{\hspace{-10pt}} | >{\hspace{-10pt}}c<{\hspace{-10pt}} | >{\hspace{-10pt}}c<{\hspace{-10pt}} | >{\hspace{-10pt}}c<{\hspace{-10pt}} | >{\hspace{-10pt}}c<{\hspace{-10pt}} | >{\hspace{-10pt}}c<{\hspace{-10pt}} | >{\hspace{-10pt}}c<{\hspace{-10pt}} | >{\hspace{-10pt}}c<{\hspace{-10pt}} | >{\hspace{-10pt}}c<{\hspace{-10pt}} |}
  \hline			
  $L_1 = L_2 = L$ & 21 & 24 & 27 & 30 & 33 & 36 & 39 & 42 & 45 \\
  \hline
  $N_{\textrm{ctrl}}$ & 49 & 64 & 81 & 100 & 121 & 144 & 169 & 196 & 225 \\
  \hline  
\end{tabular}
\caption{Domain lengths and number of controls.}\label{tab:numerics}
\end{table}\vspace{-0.25cm}
\noindent For each case in Table~\ref{tab:numerics}, we perform simulations for the unique equidistant and quasirandom (using the $2,3$-Halton sequence) actuator grids, one example of a perturbed equidistant grid, and five different random actuator arrangements, constituting $8$ grids for each of the $9$ choices of $Q$. For each simulation, we take a random initial condition with zero spatial average, containing sufficiently many unstable low modes. In contrast to the previous subsection, we allow the system to evolve without controls until it reaches the global attractor ($200$ time units suffices for the above cases), and then apply proportional point-actuated controls. For our choice of $\alpha = 150$, exponential stabilisation of the zero solution ($\overline{\eta} = 0$) is observed for all $8$ actuator arrangements with $L = 21$. At this domain size, we observe the emergence of bimodal states in the absence of controls, with the onset of chaos for slightly larger $L$.

The results of the numerical experiment are shown in Figure \ref{propcontrolsupcritscatter}. The decay rate $\lambda$ is plotted against the area predictors $|A_1 - A_1^E|$, $A_2$ and $A_3$ in panels (a--c) respectively. The equidistant grids performed the best, closely followed by the perturbed equidistant grids, both with $\lambda \approx 0.4$. The quasirandom grids, which are much more regularly spaced than the random grids, all gave exponentially decaying costs with $\lambda\approx 0.2$. The solutions controlled with random actuator grids either reached a non-trivial steady state (for which we assign $\lambda = 0$), or decayed to zero with rate no greater than $0.1$. Note that since $\kappa = 0.25$, the uncontrolled dynamics are bounded unlike in the previous subsection. Although not obvious from Figure \ref{propcontrolsupcritscatter}, we found that the proportion of random grids that resulted in exponential stability of the flat state decreased as $L$ increased, with no successes for $L=45$.

From Figure \ref{propcontrolsupcritscatter}(a), we see that the maximum attainable decay rate is a monotonically decreasing function of $|A_1 - A_1^E|$ (the points are bounded above by a monotonically decreasing curve). Actuator arrangements with $|A_1 -A_1^E| \gtrsim 10$ perform poorly, with exponential decay rates below $0.1$. For $j = 2,3$, we fit a threshold model using the data for $L=39, 42, 45$ of the form $\lambda = a_j(A_j - A_j^{\textrm{c}})^2$ for $A_j \leq A_j^{\textrm{c}}$, and $\lambda = 0$ (no exponential decay) for $A_j > A_j^{\textrm{c}}$; the constants $a_j$ and $A_j^{\textrm{c}}$ are computed using least squares fitting and optimisation. The threshold model is not appropriate for fitting with the predictor $|A_1 - A_1^E|$, since the data-points do not clearly lie on a curve, whereas it is evidently appropriate for $A_2$ and $A_3$. We obtain $\lambda = 2.1\times 10^{-3}(A_2 - 23.4)^2$ for $A_2 \leq 23.4$ and $\lambda = 1.4\times 10^{-4}(A_3 - 69.1)^2$ for $A_3 \leq 69.1$ (with $\lambda = 0$ otherwise). The least squares error may be used to quantify the ability of the areas as predictors for the success/failure of controls -- we find the least squares errors in both panels (b) and (c) to be approximately the same, thus the areas $A_2$ and $A_3$ are equally suitable predictors for the success of the controls. We expect that such measures of spacing would also be useful for the point-actuated control of other physical systems. In particular, such measures may be useful in problems where geometrical constraints are placed on the grids of actuators/observers.

For equidistant grids, exponential stability was achieved for all cases in this numerical experiment. However, such grids fail when the actuators lie at the zeros of unstable eigenfunctions (in this section, the grid spacing was finer than shortest unstable wavelength of the system). Although the results are not presented here, we found that in such situations, slightly shifting the point-actuator locations, i.e. using a perturbed equidistant actuator grid, does not prevent failure of the controls.

\subsection{Non-trivial desired states: synchronisation of chaotic dynamics.}

\begin{figure}
\begin{subfigure}{2.8in}
\caption{Costs.} 
\includegraphics[width=2.8in]{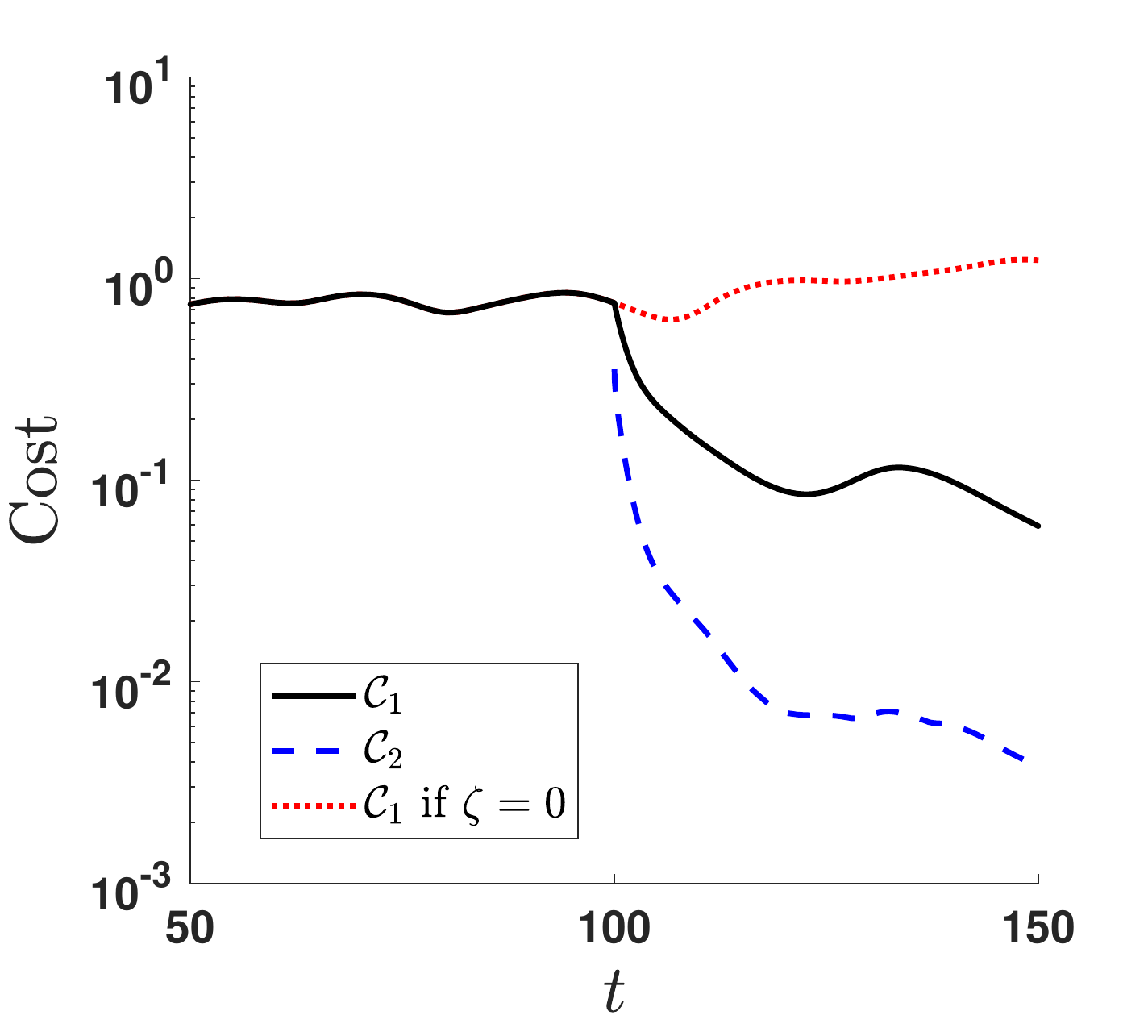}
\end{subfigure}
\begin{subfigure}{2.8in}
\caption{Projection of dynamics.}
\includegraphics[width=2.8in]{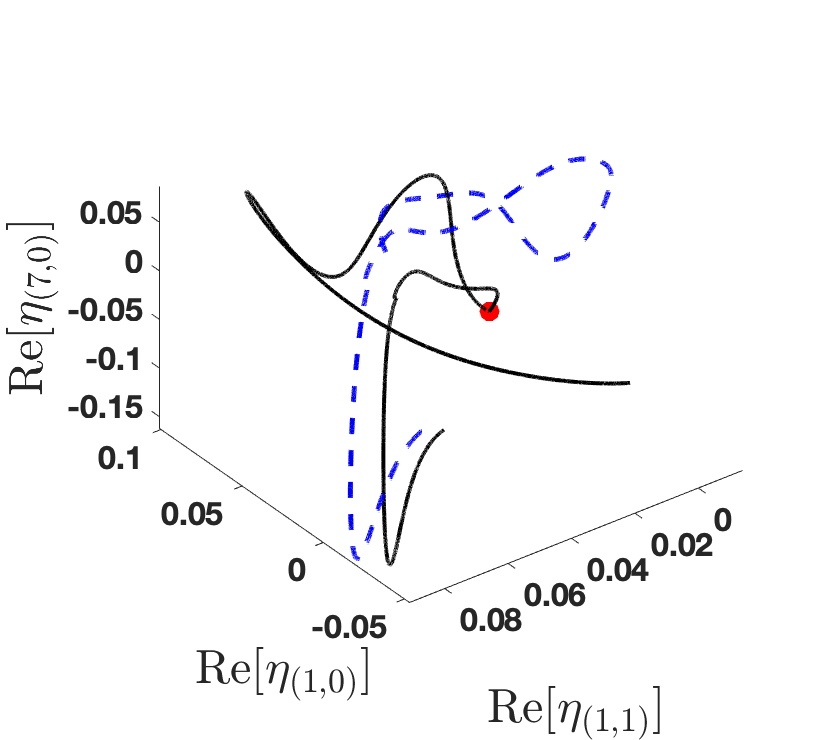}
\end{subfigure}
\begin{subfigure}{2.8in}
\caption{Costs.} 
\includegraphics[width=2.8in]{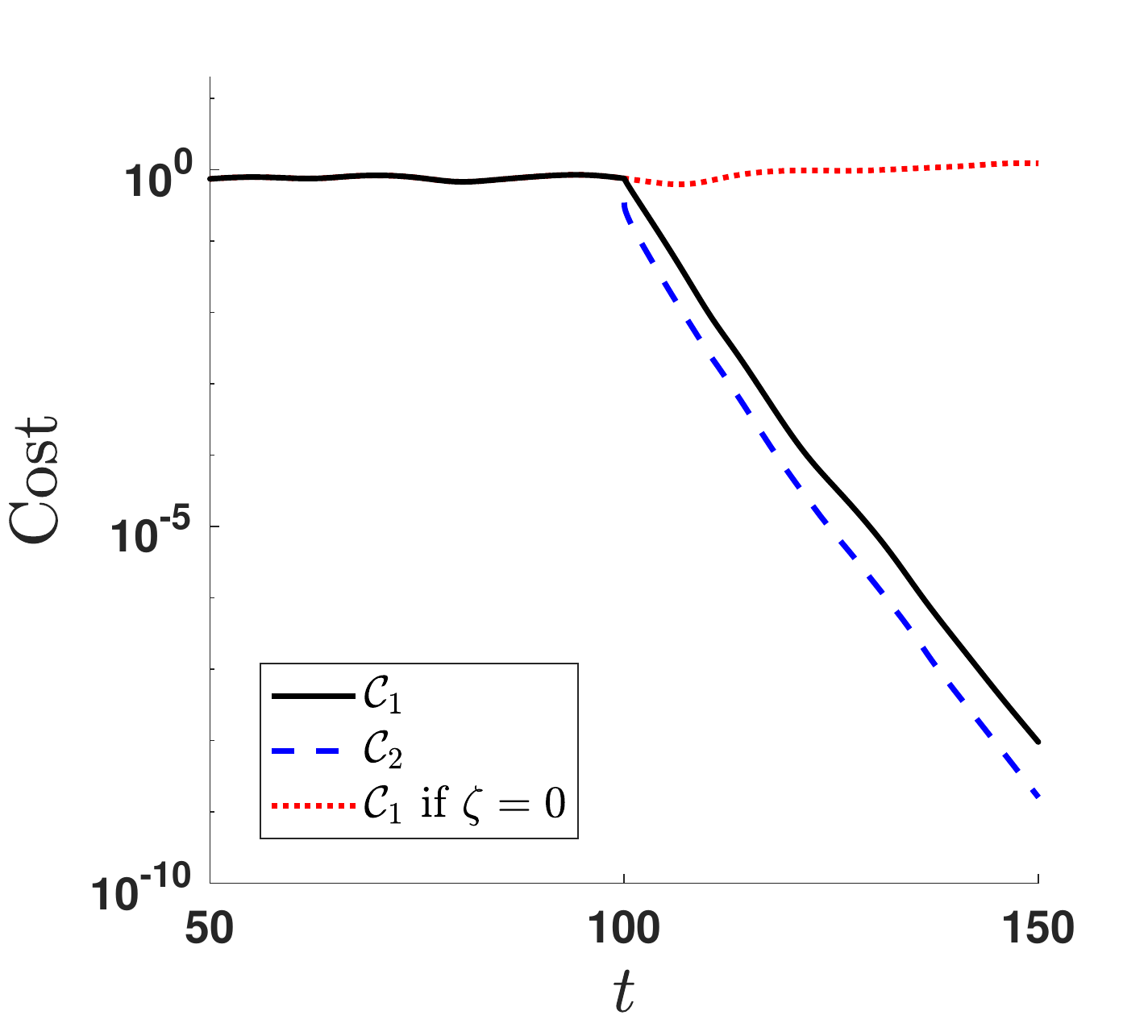}
\end{subfigure}
\begin{subfigure}{2.8in}
\caption{Projection of dynamics.}
\includegraphics[width=2.8in]{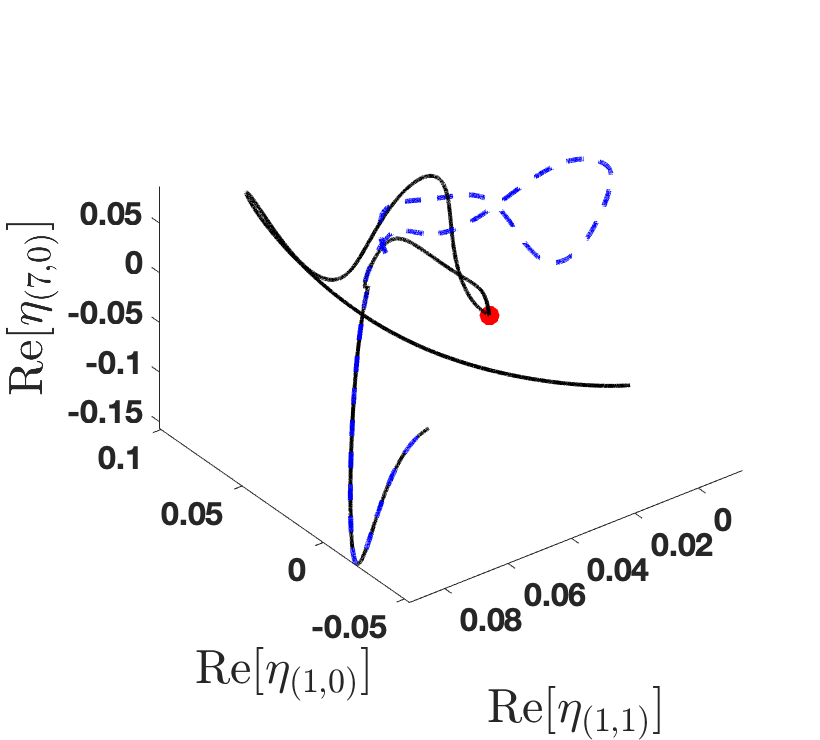}
\end{subfigure}
\caption{Synchronisation of two chaotic orbits of the 2D KSE \eqref{controlled2dks}. Panels (a,c) plot the costs for the simulations with the quasirandom and equidistant actuator grids, respectively. Panels (b,d) plot a projection of the dynamics onto the 3D space spanned by the real parts of three modes. The solid (dotted) line corresponds to the solution $\eta$ (desired state $\overline{\eta}$), and the dot denotes the time from which proportional controls are applied. The data used for this plot corresponds to $t \in [50,150]$.} \label{chaostochaos}
\end{figure}

The proportional control methodology is very robust in the sense that it allows for any choice of desired state $\overline{\eta}$, even non-solutions. For full convergence to an arbitrary non-solution, a dense set of control actuators is required, whereas if $\overline{\eta}$ is a solution of the uncontrolled equation, then it usually can be stabilised with a finite number of actuators. Non-trivial travelling waves are a popular target state for stabilisation in the KSE control literature, however, in this section, we show that even more complicated situations can be tackled, and employ proportional controls to stabilise a given chaotic orbit of the 2D KSE \eqref{controlled2dks}.

The problem of synchronising a pair of solution orbits of the 1D KSE \eqref{1dintroks} was considered by \cite{junge1999synchronization} and \cite{tasev2000synchronization}. The authors used actuators of non-zero width, and performed proportional control using local spatial averages of $\eta - \overline{\eta}$ over the actuator regions -- in the limit as the actuator/averaging width becomes small, this converges to \eqref{phiform1} with the Dirac delta actuators. They found that the number of actuators required for synchronisation scaled with the size of the periodic domain, consistent with the results of \cite{gomes2016stabilizing} and our own observations in the previous subsection. \cite{junge1999synchronization} and \cite{tasev2000synchronization} also observed that the equidistant arrangement was close to optimal -- they suggest that the discrepancy may be due to the imposed rigid boundary conditions. This study was extended to a generalised KSE by \cite{basnarkov2014generalized} with the addition of a third order dispersion term, $\gamma\eta_{xxx}$; they discussed the synchronisation of solutions to systems with different values of the dispersion strength.

In our numerical experiment, we take parameters $\kappa = 0.25$ and $L_1 = L_2 = 45$ (52 unstable modes in total) using both quasirandom and equidistant grids with $N_{\textrm{ctrl}} = 225$ actuators as in the previous subsection. We take the initial condition $\eta(\bm{x},0) = \eta_0(\bm{x})$ given in \eqref{initicond1}, and the desired state is the orbit starting from $\overline{\eta}(\bm{x},0) = 2\eta_0(\bm{x})$. Proportional controls (recall \eqref{phiform1} for the case of a non-trivial desired state) are applied with $\alpha = 5$ from time $t=100$. Figure \ref{chaostochaos} shows the results of the numerical simulations; the results for the quasirandom and equidistant arrangements are shown in panels (a,b) and (c,d), respectively. The costs for the equidistant case decay exponentially, it also appears that the costs in panel (a) decay exponentially for at least part of the controlled evolution. The projection of the orbits onto the 3D phase space with components $(\operatorname{Re}[\eta_{(1,0)}],\operatorname{Re}[\eta_{(1,1)}],\operatorname{Re}[\eta_{(7,0)}])$ tells a similar story, with the orbit of the controlled solution tracking that of $\overline{\eta}$ much more closely. See Movie 1 at \verb|https://youtu.be/kjLk9e9w5ew| for a movie of the synchronising chaotic interfaces with the quasirandom actuator grid.

\section{Feedback control with full state observations\label{SecFeedbackcontrol}}

In this section we study feedback control strategies with observation of the entire interface, which, along with the knowledge of the linearised dynamics, is employed to construct controls. This contrasts the previous section, where actuation at a point was governed by the local film thickness alone. We note that, as will be made clear later in this section, this control methodology only requires observation of a finite set of Fourier modes of the interface -- this is advantageous for the numerical implementation.

In order to apply the following theory, the assumption that $\overline{\eta}$ is an exact solution of \eqref{controlled2dks} with $\zeta = 0$ is required. Convergence to non-solutions can only be achieved transiently by such a control methodology -- see the discussion in \cite{thompson2016stabilising}. We consider the difference $w = \eta - \overline{\eta}$, which evolves according to
\begin{equation}\label{Wpdeequation112}w_t = - ww_x - (\overline{\eta}w)_x + ( - (1- \kappa)\partial_x^2 + \kappa \partial_y^2 - \Delta^2  )w + \zeta.\end{equation}
We define $A$ to be the diagonal matrix 
%indexed by $\bm{k},\; \bm{l} \in \mathbb{Z}^2$ {\color{blue} Where do we use $\bm{l}$?!}
%(the truncated matrix will be square with $(2M+1)(2N+1)$ rows and columns) 
with
%diagonal entries
\begin{equation}A_{\bm{k},\bm{k}} = (1- \kappa)\tilde{k}_1^2 - \kappa \tilde{k}_2^2 - | \bm{\tilde{k}}|^4,\end{equation}
and $B$ denotes the control actuator matrix with entries $B_{\bm{k},j} = b_{\bm{k}}^j$. We assume that the controls $\phi^{j}$ depend linearly on the Fourier coefficients of $w$ through a constant matrix $K$ (to be chosen) by the relation
\begin{equation}\label{feedbackfullphij}\phi^j = \sum_{\bm{l}\in\mathbb{Z}^2} K_{j,\bm{l}} w_{\bm{l}} = [K \bm{w}]_j,\end{equation}
where $\bm{w}$ is the vector containing the Fourier coefficients of $w$.
% and the superscript $T$ denotes the transpose. 
 This is known as a linear state feedback control law, and $K$ is known as the feedback gain matrix.
%Again, we take $K$ to be a two-dimensional matrix with the same dimensions as $B$ rather than a three-dimensional one. 
The control $\zeta$ has Fourier coefficients
\begin{equation}\zeta_{\bm{k}} = \sum_{j = 1}^{N_{\textrm{ctrl}}} \sum_{\bm{l}\in\mathbb{Z}^2} B_{\bm{k},j}   K_{j,\bm{l}} w_{\bm{l}} = [BK \bm{w}]_{\bm{k}}.\end{equation}
With this, we may rewrite (\ref{Wpdeequation112}) in terms of matrices as
\begin{equation}\label{Wequation111}\bm{w}_t = \bm{\nu} + (J + A + BK) \bm{w},\end{equation}
where $\bm{\nu}$ is the vector with entries being the Fourier coefficients of $- ww_x$, and $J$ is the matrix such that $J\bm{w}$ gives the Fourier coefficients of $- (\overline{\eta}w)_x$. The entries of $J$ may be computed as $J_{\bm{k},\bm{l}} = - i \tilde{k}_1 \overline{\eta}_{\bm{k} - \bm{l}}$. The task now is to choose the feedback gain matrix $K$ so that this nonlinear infinite-dimensional dynamical system is stable about $\bm{w} = \bm{0}$. It is well known that the characterisations of stability of linear infinite-dimensional dynamical systems are more intricate than the finite-dimensional case \citep{zabczyk2009mathematical}, and further complications arise when nonlinear terms are added. We now describe two methodologies that are applicable to this problem, having both been considered for the 1D case.

\cite{doi:10.1137/140993417} showed for the corresponding 1D problem that the full infinite-dimensional nonlinear system may be controlled by ensuring the stability of a truncation of the linearised system. Linearising the problem about the desired state $\overline{\eta}$, i.e. linearising \eqref{Wequation111} about $\bm{w} = \bm{0}$ (or \eqref{Wpdeequation112} about $w=0$), gives 
\eqref{Wequation111} with $\bm{\nu} = \bm{0}$.
%The linearisation of (\ref{Wpdeequation112}) about $w = 0$ (equivalent to the linearisation of the original equation \eqref{controlled2dks} for $\eta$ about $\overline{\eta}$) 
%is
%\begin{equation}\label{Wpdeequation113linear}w_t = - (\overline{\eta}w)_x + ( - (1- \kappa)\partial_x^2 + \kappa \partial_y^2 - \Delta^2  )w + \zeta.\end{equation}
%can be written as the matrix system
%\begin{equation}\label{Wequation114linear}\bm{w}_t = (J + A + BK) \bm{w},\end{equation}
%(the contribution from $\bm{\nu}$ vanishes from the linearised dynamics). 
We let $A^n$ denote the truncation of the matrix operator $A$ to the modes with $|k_1|,|k_2| \leq n$, and we have similar notations for the other matrices and vectors. The analysis of \cite{doi:10.1137/140993417} can be lifted to our 2D setting, and their Theorem 5.1 may be recast as:
\begin{thm}\label{MorrisThm1}(Theorem $5.1$ in \cite{doi:10.1137/140993417})
Consider the sequence of approximations of the linear matrix problem \eqref{Wequation111} with $\bm{\nu} = \bm{0}$ defined by the Galerkin truncation in Fourier space onto the modes with $|k_1|,|k_2| \leq n$,
\begin{equation}\label{eq:MorrisThm}
\bm{w}^{n}_t = \left(J^{n} + A^{n} + B^{n} K^{n} \right) \bm{w}^{n}.
\end{equation}
Assume that there exists a convergent sequence of matrices $K^{n}$ which stabilise the problem~\eqref{eq:MorrisThm}, and such that the limit $K$ exponentially stabilises \eqref{Wequation111} with $\bm{\nu} = \bm{0}$. Then for sufficiently large ${n}$, the controller $K^{n}$ stabilises the full nonlinear problem \eqref{Wequation111}.
\end{thm}
Thus, to apply this numerically to our problem, we take a sufficiently large truncation $n$ and compute a (possibly time varying) matrix $K^{n}$ such that \eqref{eq:MorrisThm} is stable -- this is then extended with zeros to obtain a viable choice of $K$. We expand upon how such a $K^{n}$ is chosen later, but we note that it must possess certain symmetries to ensure that the numerical solution remains real-valued. This method is computationally feasible for a restricted choice of desired states. If $\overline{\eta}$ is a steady state, then the matrix $J^{n}$ is constant in time ($J^{n}$ is zero if $\overline{\eta} = 0$), and so $K^{n}$ is constant in time and only needs to be computed once, giving a static feedback law. If $\overline{\eta}$ is an exact travelling wave solution with period $\tau$, then a value of $K^{n}$ needs to computed for all time-steps in $[0,\tau]$ as $J^{n}$ varies (the matrices $K^{n}$ should also vary continuously in time); this results in a dynamic feedback law. For choices of $\overline{\eta}$ with even more complicated dynamics, such as quasi-periodicity or chaos, $K^{n}$ needs to be computed at each time-step for the entire time interval for which controls are applied; this is computationally excessive and unfeasible in our case.

Another method was introduced by \cite{gomes2016stabilizing} in their study of the 1D problem. The authors linearise about the zero state regardless of the desired state $\overline{\eta}$ (this methodology coincides with the former one for $\overline{\eta} = 0$). The matrix $K^{n}$ is chosen to stabilise the linear system
\begin{equation}\label{Wequation1152linear}\bm{w}_t = \left( A^{n} + B^{n} K^{n} \right) \bm{w},\end{equation}
under a further restriction which averts the problem caused by the nonlinearity, as described next. 
For $\lambda > 0$, we define $\lambda$-stability for a complex matrix $C$ to be the property that for any complex vector $\bm{v}$,
\begin{equation}\label{lambdastability1}(\bm{v}^*)^T C \bm{v}  \leq - \lambda  (\bm{v}^*)^T \bm{v} = - \lambda \| \bm{v} \|_{\ell^2},\end{equation}
where the $*$ denotes taking the complex conjugate. If $C$ is a normal matrix, satisfying $C^{*}C^{T} = C^{T} C^{*}$, condition \eqref{lambdastability1} is equivalent to the real parts of all eigenvalues of $C$ being bounded above by $-\lambda$. 
Following the same argument as \cite{gomes2016stabilizing}, in addition to stabilising \eqref{Wequation1152linear}, the truncated feedback gain matrix $K^{n}$ must be chosen such that the bracketed term in \eqref{Wequation1152linear} is $\lambda$-stable
for 
\begin{equation}\label{lambdastablebound1}\lambda + \frac{1}{2} \inf_{\bm{x} \in Q} \overline{\eta}_x > 0.\end{equation}
Then, if the truncation is suitably large, the full nonlinear system \eqref{Wequation111} is exponentially stabilised with decay rate given by the left hand side of \eqref{lambdastablebound1}. In practice, we replace the $\lambda$-stability requirement with the weaker condition of having the real parts of all eigenvalues bounded above by $-\lambda$, as proved successful in \cite{gomes2016stabilizing} -- these are not equivalent conditions for our problem since the bracketed term in \eqref{Wequation1152linear} is non-normal. We find that the eigenvalue condition alone is sufficient to stabilise the systems we consider numerically, and $\lambda$-stability may be too strong a requirement.

Computing a feedback gain matrix $K^{n}$ to yield desired eigenvalues for the linear systems (\ref{eq:MorrisThm},\ref{Wequation1152linear}) is a pole placement problem. For this we use the \textsc{Matlab} function \textit{place}. This procedure is difficult to carry out in the complex Fourier mode framework, since in this basis, the entries of $K^{n}$ are complex yet must satisfy symmetry requirements so that the resulting forcing is real-valued. Thus, we translate the problem into real-valued, trigonometric basis functions by applying linear transformations to the matrices $J^{n}$, $A^{n}$, and $B^{n}$. Given the matrices $J^{n}$, $A^{n}$, $B^{n}$ (in the new basis), and a vector of desired eigenvalues, \textit{place} computes a suitable real matrix $K^{n}$ such that the linearised system possesses the desired eigenvalues. This procedure is robust to changes in $A^{n}$ and $B^{n}$ \citep{kautsky1985robust}. The matrix $K^{n}$ may then be transformed into the original complex Fourier mode basis and used in numerical simulations. The precise details of how the desired eigenvalues are chosen is given in the following text, yet in broad terms, we choose the same set of eigenvalues as in the uncontrolled system, but replacing those above a given threshold with a negative value. It is important to note that, with this procedure of choosing the desired spectrum, the matrices $K^n$ do not form a convergent sequence as described in Theorem \ref{MorrisThm1}. For large values of $n$, the method breaks down, with $K^n$ having large entries, resulting in a heavily sensitive problem. Theorem \ref{MorrisThm1} makes no restriction on the exact form of the spectra, so a more appropriate algorithm for this situation would be one that places importance on having all eigenvalues below a threshold, while ensuring that the entries of the feedback gain matrix are relatively small. If we were to employ the methodology of \cite{doi:10.1137/140993417} for more complicated choices of $\overline{\eta}$, we would need a pole-placement algorithm which computes the feedback gain matrix rapidly, and ensures that $K^n$ varies continuously as $J^n$ varies. However, such extensions to the algorithm are beyond the scope of the current paper.

We now present numerical experiments testing the effectiveness and applicability of such controls in three different situations.

\subsection{Controlling to the trivial flat solution\label{feedbackflatsubsec}}

\begin{figure}
\begin{subfigure}{2.8in}
\caption{Costs.} 
\includegraphics[width=2.8in]{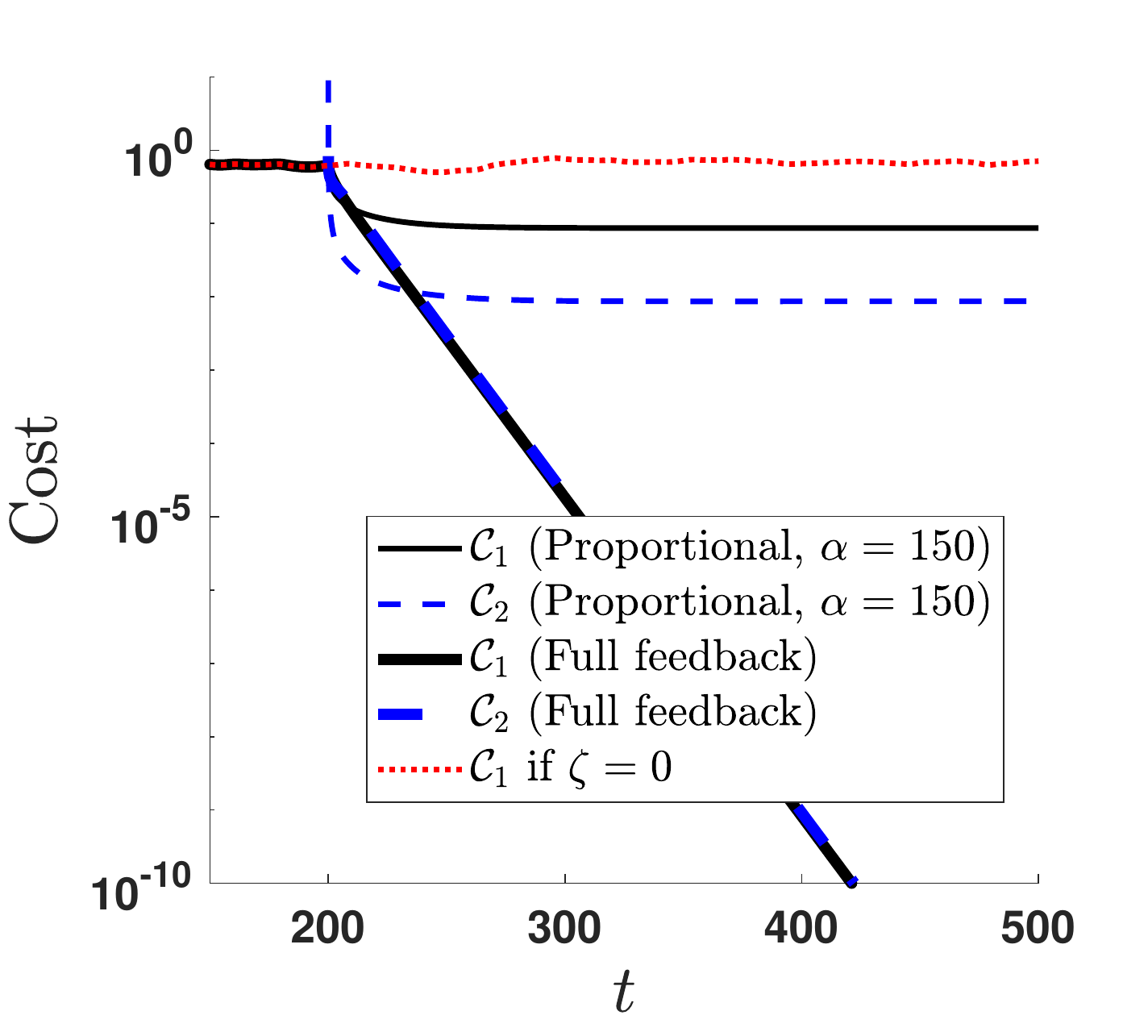}
\end{subfigure}
\begin{subfigure}{2.8in}
\caption{Actuation strength.}
\includegraphics[width=2.8in]{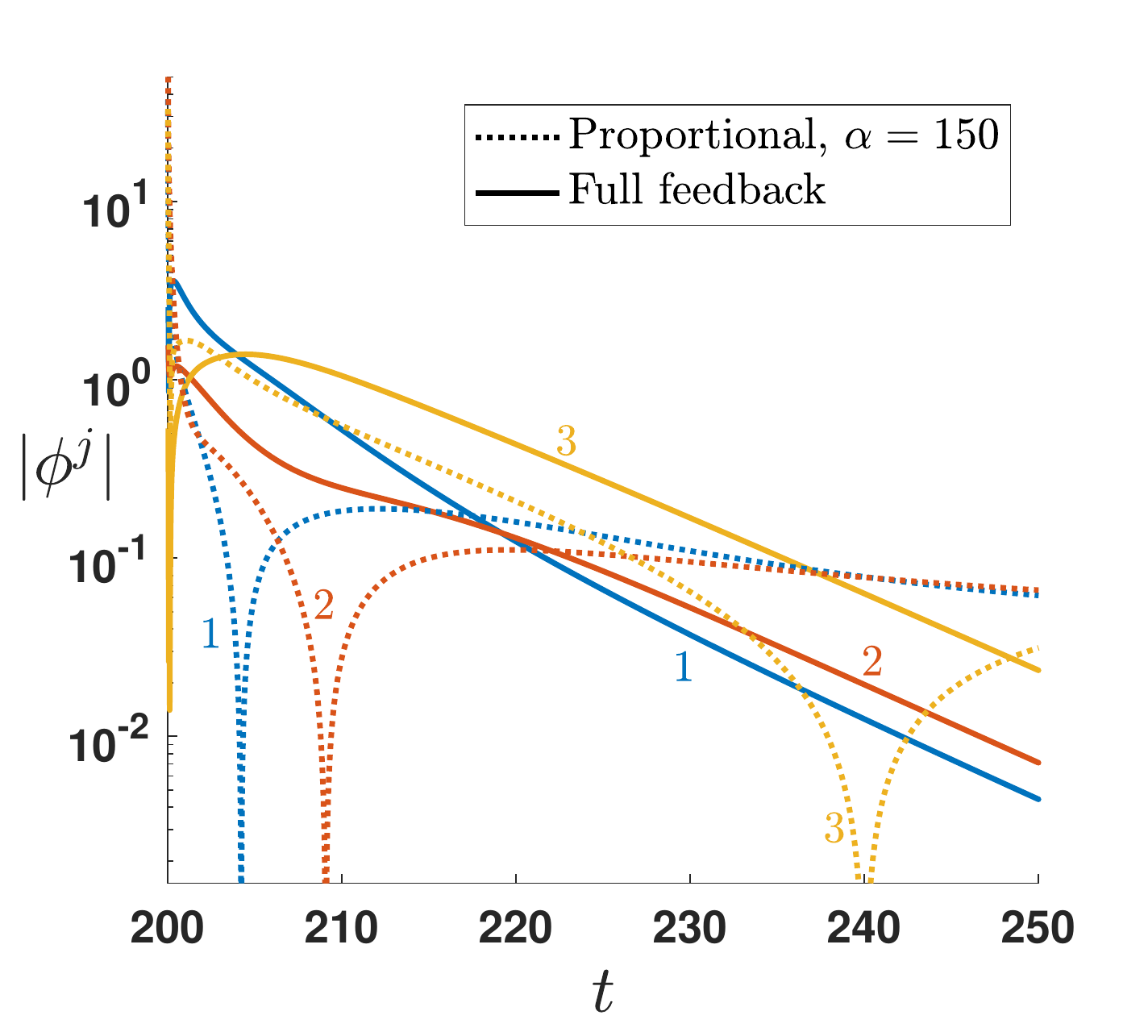}
\end{subfigure}
\caption{Exponential stabilisation of zero solution with full feedback controls. Panel (a) shows the evolution of the costs $\mathcal{C}_1$ and $\mathcal{C}_2$; the thick (thin) lines correspond to the results for full feedback (proportional) controls with $\alpha = 150$. The uncontrolled solution cost is also plotted. Panel (b) compares the control strength for the same three actuators across the proportional and full feedback control strategies.} \label{feedbackcostssupcrit}
\end{figure}
We first apply feedback controls to stabilise the flat film solution, taking $\overline{\eta} = 0$. In this case, the methodologies in \cite{doi:10.1137/140993417} and \cite{gomes2016stabilizing} agree. We take parameters $\kappa = 0.25$, $L_1 = L_2 = 42$, $N_{\textrm{ctrl}} = 196$, and a random actuator grid used in the simulations in subsection \ref{CompofGridssubsec}. With proportional controls and $\alpha = 150$, exponential stability of the zero solution was not achieved for this particular actuator arrangement. As in subsection \ref{CompofGridssubsec}, controls are applied after $200$ time units, and the (random field) initial condition is the same as used there also. 
%This particular set-up is used for all the control strategies considered in this work as a point of comparison. 
Using a truncation of $n = 19$, which importantly covers the linearly unstable modes, the matrix $K^{n}$ is computed using the \textsc{Matlab} function \textit{place} so that the eigenvalues of the linearised system \eqref{Wequation1152linear} are at most $-0.1$. More precisely, the eigenvalues of the linearised system are unchanged if they are less than $-0.1$, with the rest replaced by $-0.1$. This decay rate improves upon any of the decay rates obtained by random grids and proportional controls in subsection \ref{CompofGridssubsec}.

The evolution of the costs are shown in Figure \ref{feedbackcostssupcrit}(a), and it can be seen that exponential stabilisation is achieved with decay rate $0.1$ (correct to 7dp) for the full feedback controls; this is because $-0.1$ was chosen to be the largest eigenvalue of the truncated linearised system. We find that $\mathcal{C}_2(200)$ (the initial control cost) is an order of magnitude smaller for full feedback control than for proportional control. We found that $\mathcal{C}_1 \approx \mathcal{C}_2$ for the feedback controlled case. In Figure \ref{feedbackcostssupcrit}(b) we plot the actuation strengths $|\phi^j|$ for the same three actuators across the proportional and full feedback control strategies. The actuation strengths for the full feedback control case (shown with solid lines) all decay exponentially after an initial transient phase, but the proportional controls ``over-control" the interface since the controls change sign (indicated by the singularities in the $\log$-linear plot). 

The full feedback controls are thus shown to perform well for the case of $\overline{\eta} = 0$, even with actuator arrangements that are far from optimal. The costs decayed exponentially with the predicted value, and the solution converged to the desired state to machine precision. In the next two subsections, we see how choosing a non-trivial desired state complicates the problem.

\subsection{Controlling to non-trivial steady states}

\begin{figure}
\centering
\begin{subfigure}{2.8in}
\caption{Streamwise slice of steady state.}
\includegraphics[width=2.8in]{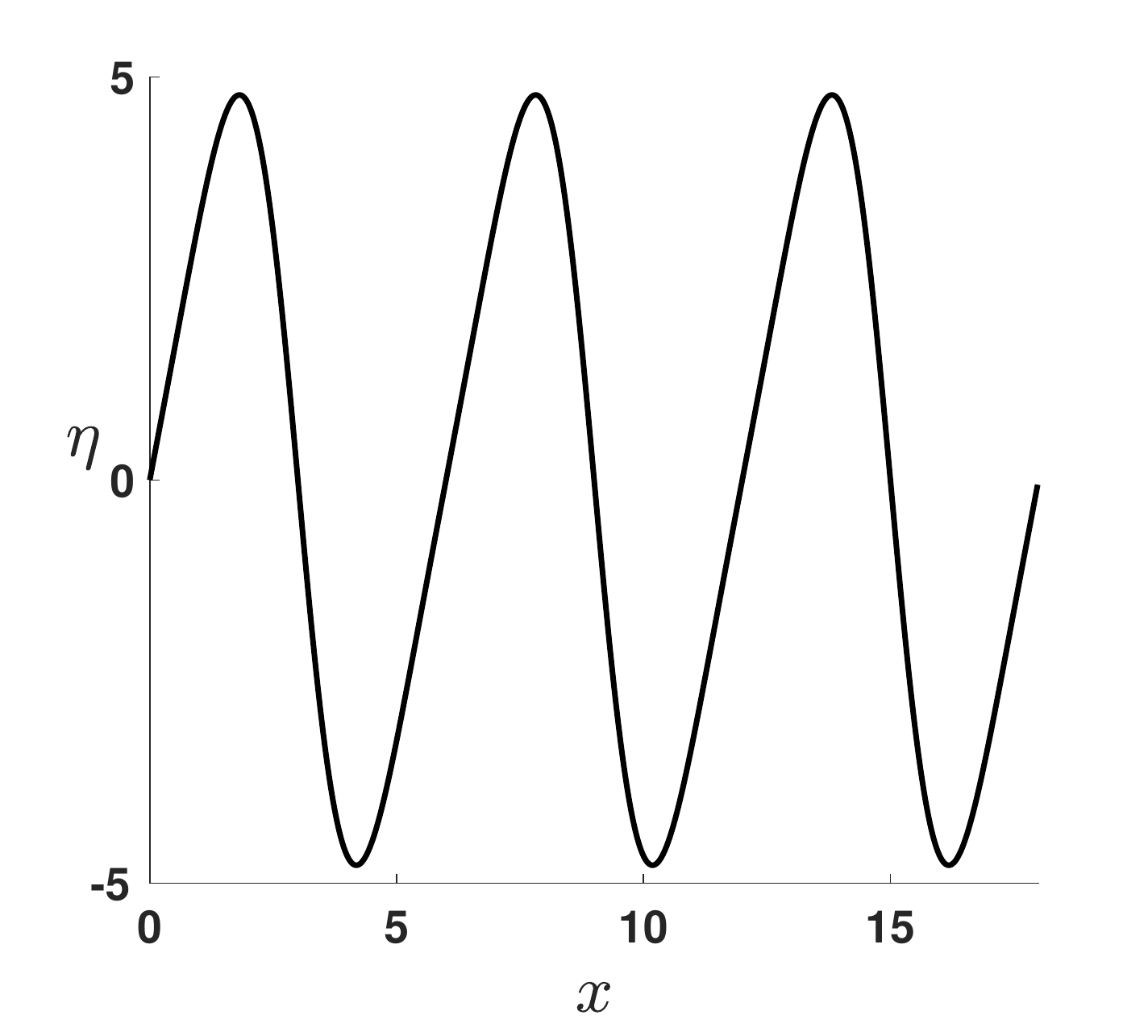}
\end{subfigure}
\begin{subfigure}{2.8in}
\caption{Costs.}
\includegraphics[width=2.8in]{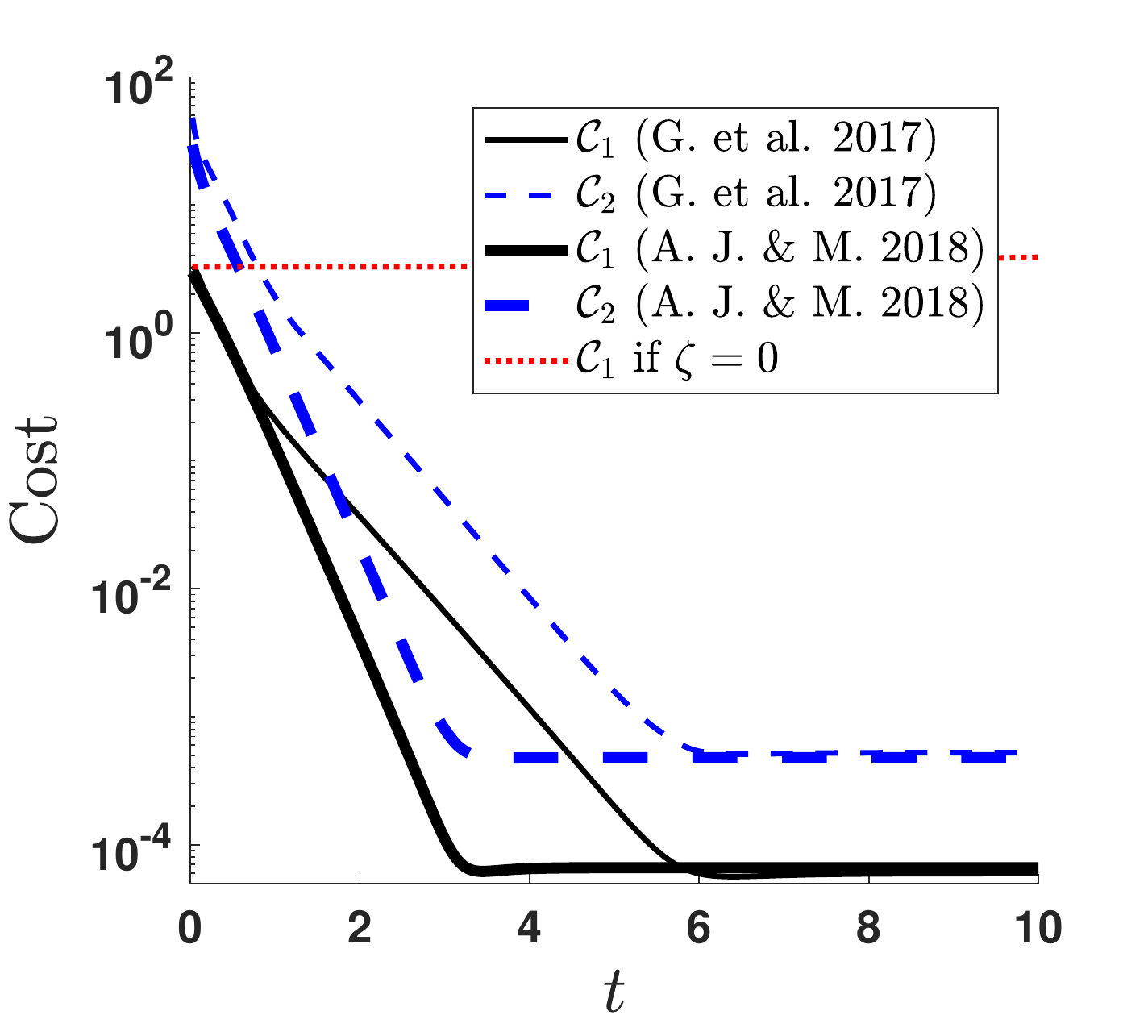}
\end{subfigure}
\caption{Stabilisation of an exact steady state solution (details of parameters are given in the text). Panel (a) shows the profile of the one-dimensional steady state. Panel (b) plots the costs $\mathcal{C}_{1}$ and $\mathcal{C}_{2}$ for the applied full feedback controls for both the methodologies of \cite{doi:10.1137/140993417} and \cite{gomes2016stabilizing}, with $\mathcal{C}_1$ for $\zeta = 0$ included for reference.} \label{sswavefig}
\end{figure}

In this subsection, we consider the stabilisation of a non-trivial steady solution of \eqref{controlled2dks} with full feedback controls; the methodologies of \cite{doi:10.1137/140993417} and \cite{gomes2016stabilizing} differ in this case, yet they are both computationally feasible. For this numerical experiment, we take the parameters $\kappa = -0.5$, $L_1 = L_2 = 18$, and use $N_{\textrm{ctrl}} = 100$ randomly located actuators (the arrangement used for numerical experiments in subsection \ref{subsectionsuppressgrohanging}, see Figure \ref{hangingpropcostsC1C2}). We choose $\overline{\eta}$ to be an unstable 1D steady state which is constant in $y$; the profile is plotted in Figure \ref{sswavefig}(a). Such 1D steady states may be computed with ease using the continuation and bifurcation software \textsc{AUTO-07P}, for example. This exact solution of the uncontrolled 2D KSE is unstable to transverse perturbations, and is even unstable to streamwise perturbations; chaos is already prevalent in the 1D KSE with these parameters. Recall that, for the application of the method of \cite{gomes2016stabilizing}, we must satisfy the $\lambda$-stability condition \eqref{lambdastablebound1}. For the steady state shown in Figure \ref{sswavefig}(a), the infimum of the $x$-derivative can be computed to be approximately $-6.762$. Assuming the equivalent condition for normal matrices, we satisfy the $\lambda$-stability condition by choosing the eigenvalues of the linearised system \eqref{Wequation1152linear} to be bounded above by $-3.5$. Since the number of repeated eigenvalues cannot be more than the number of controls (the pole placement algorithm fails if the multiplicity of a desired eigenvalue is greater than the rank of the matrix $B^n$), the desired spectrum for the controlled problem is chosen similarly as in the previous subsection, but with the eigenvalues larger than $-3.5$ replaced with $-(3.5 + 0.1U)$ where $U \sim \textit{Unif}(0,1)$. To ensure a robust comparison, we chose the eigenvalues for our computations using the method of \cite{doi:10.1137/140993417} in a similar way, replacing the real parts of the complex eigenvalues of $J^n + A^n$ (note that the eigenvalues of $A^n$ alone are real-valued) which exceed $-3.5$ with values below this bound. We use a mode truncation of $n = 9$ for the construction of the feedback gain matrix. The costs obtained are plotted in Figure \ref{sswavefig}(b). The costs corresponding to the methodology of \cite{doi:10.1137/140993417}, shown with thick lines in the figure, decay at the expected rate of $3.5$ approximately. For the method of \cite{gomes2016stabilizing} shown with thin lines in Figure \ref{sswavefig}(b), the costs initially proceed as in the previous case, before switching to a less extreme decay rate at around $t = 1$. In both cases, $\mathcal{C}_1$ decays to a plateau at approximately $6\times 10^{-5}$ (with a corresponding plateau for $\mathcal{C}_2$). This behaviour does not appear to be due to inaccuracy of the computed steady state ($\overline{\eta}$ was computed to much higher accuracy), nor is it an effect of violation of the $\lambda$-stability criterion since it is seen for both methodologies. We also confirmed that it was not due to numerical error (time or space discretisation). Thus, we attribute it to the effect of nonlinearity and the truncation of the controlled system.

%The inaccuracy appears in $k_1 = 3, 6$ modes...
%
%
%{\color{red} Could be the difference in norm between two adjacent time step solutions? 
%
%
%}

\subsection{Controlling to travelling waves}

\begin{figure}
\centering
\begin{subfigure}{2.8in}
\caption{Streamwise slice of travelling wave.}
\includegraphics[width=2.8in]{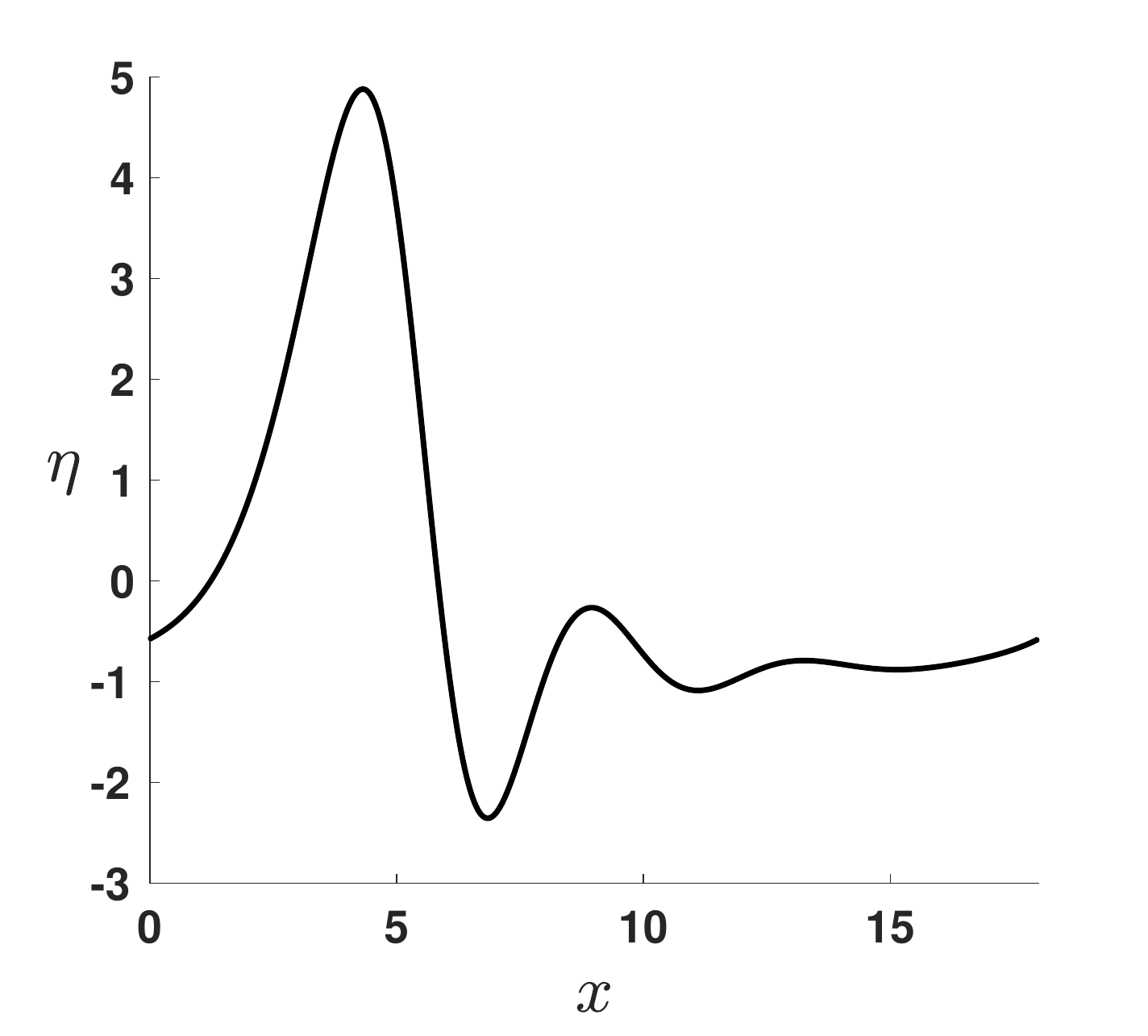}
\end{subfigure}
\begin{subfigure}{2.8in}
\caption{Costs.}
\includegraphics[width=2.8in]{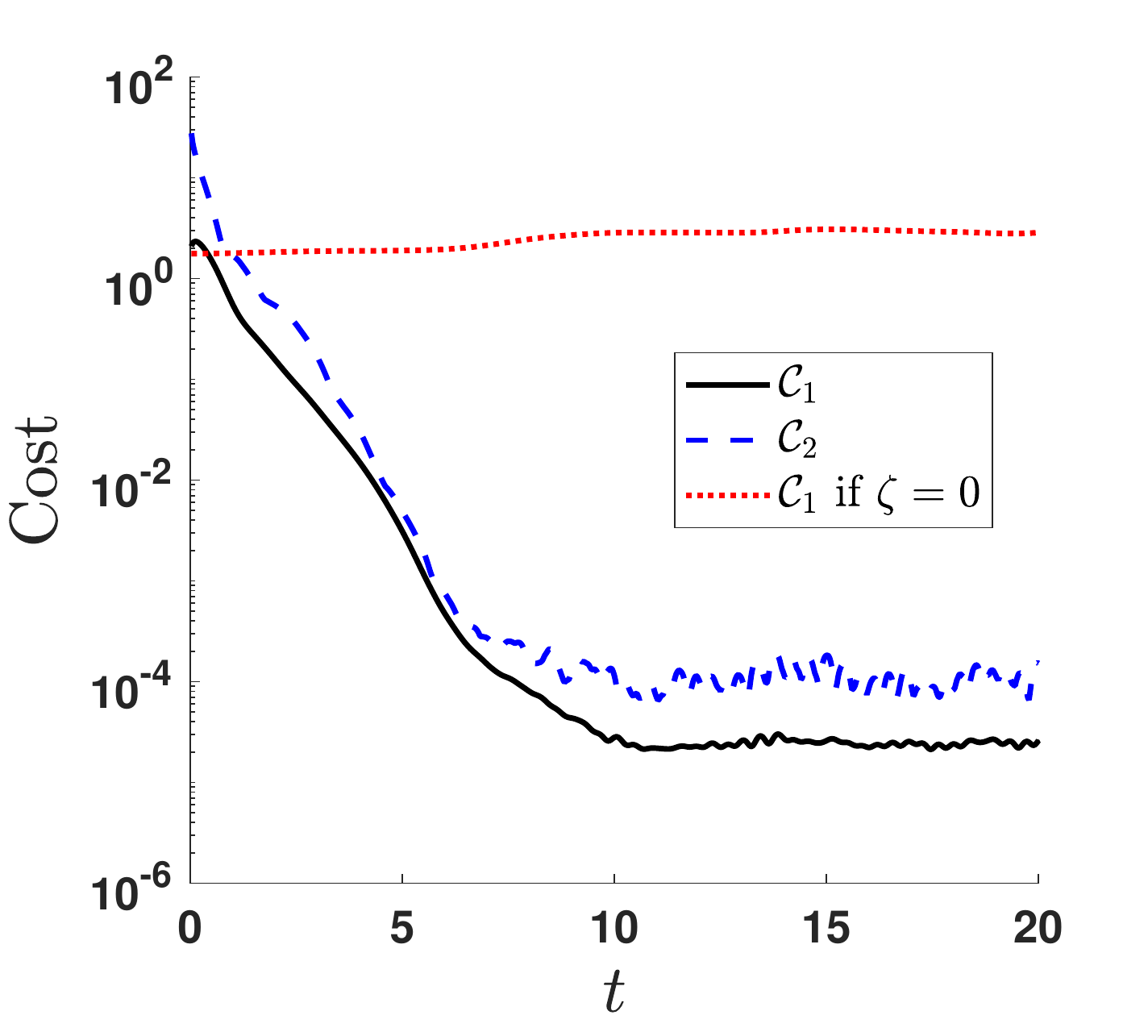}
\end{subfigure}
\begin{subfigure}{2.8in}
\caption{Controlled solution at $t=0.75$.}
\includegraphics[width=2.8in]{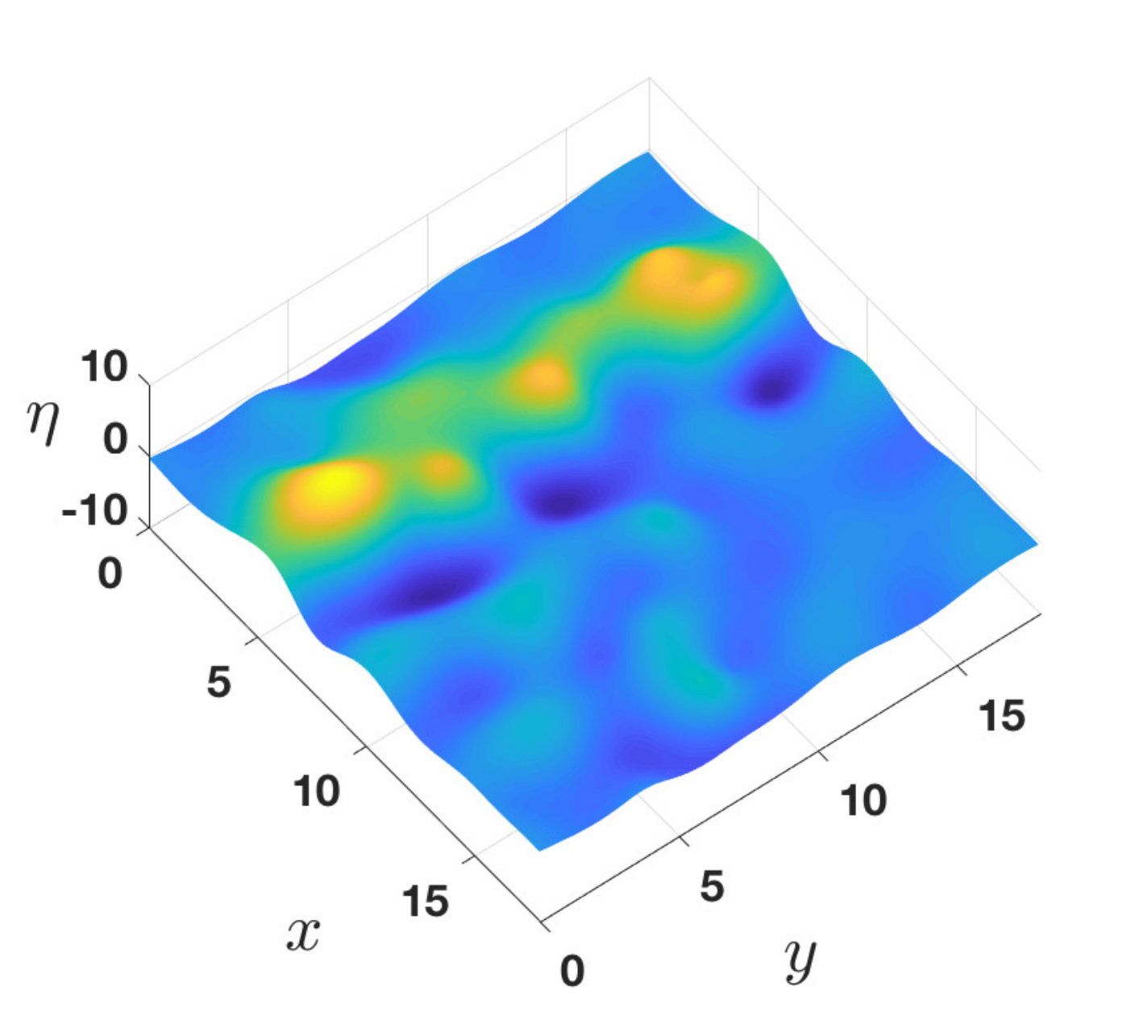}
\end{subfigure}
\begin{subfigure}{2.8in}
\caption{Controlled solution at $t=2.5$.}
\includegraphics[width=2.8in]{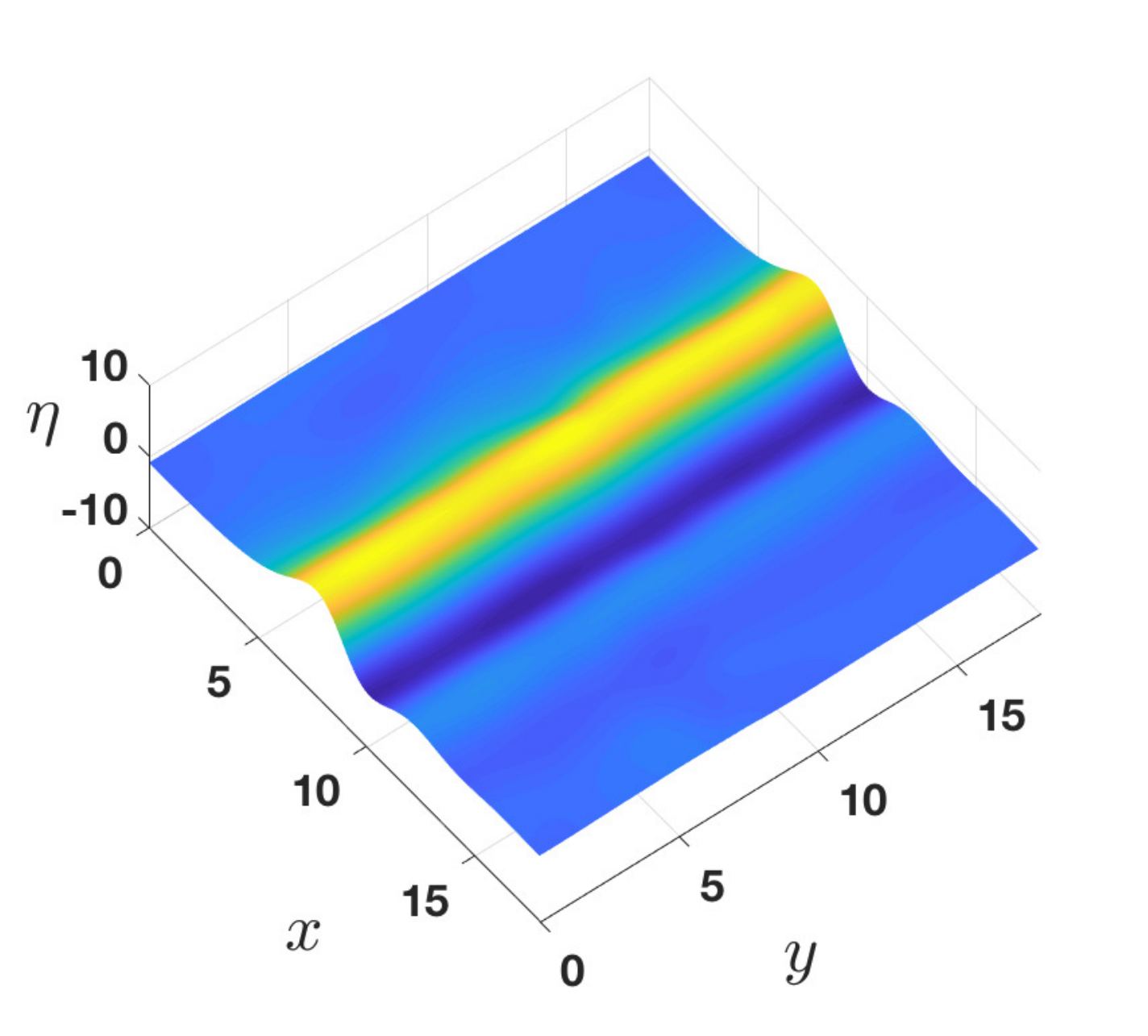}
\end{subfigure}
\caption{Stabilisation of an exact travelling wave solution (details of parameters are given in the text). Panel (a) shows the profile of the one-dimensional travelling wave which moves in the streamwise direction with speed $c \approx 1.356$ (downstream). Panel (b) plots the costs $\mathcal{C}_{1}$ and $\mathcal{C}_{2}$ for the applied full feedback controls, and gives $\mathcal{C}_1$ in the uncontrolled case for reference. Panels (c,d) show the controlled interface at times $t=0.75, 2.5$, respectively.} \label{travwavefig}
\end{figure}

For this we take the same parameters as in the previous subsection, but use a random actuator grid with $N_{\textrm{ctrl}} = 50$ (an arrangement used in subsection \ref{subsectionsuppressgrohanging}). We also consider the initial condition \eqref{initicond1}. A cross-section of the desired 1D travelling wave is shown in Figure \ref{travwavefig}(a) (as in the previous subsection, this is unstable to both streamwise and transverse perturbations). The methodologies of \cite{gomes2016stabilizing} and \cite{doi:10.1137/140993417} differ in this case of non-trivial $\overline{\eta}$, and we performed numerical simulations for the former case here since it is simpler to realise numerically. The infimum of the $x$-derivative of the desired wave shown in Figure \ref{travwavefig}(a) can be computed to be $-4.802$ approximately, thus in the aim of satisfying the $\lambda$-stability bound \eqref{lambdastablebound1}, we ensure that the eigenvalues of the linearised system \eqref{Wequation1152linear} are at most $-2.5$. We do this in a similar way to the previous numerical experiment, choosing $K^{n}$ so that we retain the same spectrum as in the uncontrolled problem, but with eigenvalues larger than $-2.5$ replaced with $-(2.5 + 0.1U)$ where $U \sim \textit{Unif}(0,1)$. The additional random variable term removes the possibility of repeated eigenvalues, which we found to be problematic for \textit{place} with this particular choice of parameters. 
For the results we present here, the largest eigenvalue was approximately $-2.503$.

We tested mode truncations of $n = 14$, $19$ and $24$. The latter choice yielded a feedback gain matrix with very large entries, and was not useful for control; the matrices $K^n$ with such large truncations are numerically expensive to compute, and may be unusable due to our method of choosing eigenvalues for the linearised system (as discussed previously). The choices of $n=14,19$ performed well, with the results for the latter case plotted in Figure \ref{travwavefig}(b) -- it can be seen that the costs plateau after an initial phase of exponential decay. Figure \ref{travwavefig}(c,d) show snapshots of the controlled interface at $t=0.75$ and $t=2.5$, respectively. The evolution of the controlled interface is shown in Movie 2 found at \verb|https://youtu.be/uEYoiD_klpk|.

Full feedback controls performed well in all cases, where a large number of modes are unstable and chaos or exponential growth is prevalent in the uncontrolled dynamics. The truncated linear system from which the feedback gain matrix is computed should include the linearly unstable modes; we were able to do this for all our examples. We remark that the method will break down with too many unstable modes (requiring a large $n$), as the problem becomes too high-dimensional for \textit{place}. A different pole placement methodology would be useful for these cases. The authors also considered controls with dynamical observers as was utilised for the 1D KSE~\eqref{1dintroks} in \cite{christofides1998feedback,armaou2000feedback}, and for the control of a 1D Benney equation by \cite{thompson2016stabilising}. We found the method to be largely unsuccessful, and that the results did not improve on those provided by proportional controls.

\section{Conclusions\label{ConcSec}}

In this paper we considered the feedback control of a multidimensional Kuramoto--Sivashinsky equation using point-actuators. The equation yields steady, travelling, time-periodic and quasi-time-periodic waves, as well as chaotic and unbounded solutions depending on the parameter regime. We applied two closed-loop control strategies, proportional control and feedback control with full state observations. For proportional control, we investigated the limitations of the method depending on the strength, number, and placement of the actuators/observers. The controls were able to prevent the unbounded growth of the interface, and exponentially stabilise a desired state. We used three measures of the actuator spacings for the grids, and found that they were strongly correlated with the decay of the controlled system, taking a maximum decay rate for equally spaced actuator arrangements. The proportional controls performed well for non-trivial desired states; we were able to synchronise two chaotic orbits of the system. We note that knowledge of the governing dynamics is not necessary to apply proportional controls, and thus may be applicable in experiments. For this purpose, it would be of interest to investigate the use of phase-shifted controls as done by \cite{thompson2016stabilising} for long-wave thin film models in the 1D setting -- although not appearing in our model, there will often be a space and time-lag involved in the control of fluid systems.

We used feedback control with full state observations to stabilise the zero solution, a steady state, and travelling wave solution. Feedback gain matrices were constructed following the methodologies of \cite{doi:10.1137/140993417} and \cite{gomes2016stabilizing} which were considered for the 1D problem, the former being an analytical result and the latter being more easily implemented numerically for non-trivial desired states; we found that both methods performed well. Although the full interface is observable, only a finite-dimensional subset of the Fourier modes is required to construct the feedback gain matrix. Furthermore, knowledge of the governing equation must be known \textit{a priori}. Current work by the authors involves data-driven control strategies for both deterministic and stochastic evolution equations; the aim is to achieve similar levels of success as feedback control with full state observations, without knowledge of the governing equations and limited observability.

The success of controls constructed for interfacial evolution equations when applied to more complicated models (for the same system) is unknown. In the thin film scenario, the authors are investigating the possibility of controlling the interface for a flow governed by the full Navier--Stokes equations using controls constructed for long-wave models. Even if the interface is successfully controlled, how will the controls affect the flow in the bulk? Interface waviness is useful in heat transfer applications, partially due to recirculation regions located in the wave crests. If controls are applied to drive the interface of a subcritical Reynolds number flow to such a wavy desired state, it is not necessarily clear that the bulk flow will recover the flow recirculation regions.

%{\color{red}
%
%
%
%\begin{itemize}
%\item Cerpa -- non-periodic problems, boundary control
%\item data assimilation 
%\item Determining modes/nodes.
%\item Non-normality (Stephen paper)
%\item Noisy observations
%\end{itemize}
%
%
%}

\section*{Acknowledgments}

R.J.T. gratefully acknowledges a PhD studentship from EPSRC. S.N.G. is supported by the Leverhulme trust via the Early Career Fellowship ECF-2018-056, and EPSRC grants EP/K034154/1 and EP/L020564/1. The authors would also like to acknowledge Professor Demetrios T. Papageorgiou and Professor Grigorios A. Pavliotis for their helpful discussions.

\bibliographystyle{imamat}
\bibliography{Controlbib}

\end{document}